\documentclass[12pt]{amsart}
\usepackage{amsfonts,amssymb,amscd,amsmath,epsfig}
\usepackage{latexsym}
\usepackage{mathrsfs}
\usepackage{tocvsec2}

\usepackage{titlesec}
\titleformat{\section}[runin]{}{\bf\thesection.}{0.5em}{}[.]

\usepackage[dvipdfm,
hyperindex=true,pagebackref=true,bookmarks=true,
colorlinks=true,linkcolor=blue,citecolor=red]{hyperref}
\def\~{{\rm --}}

\newcommand{\comment}[1]{}

\title [Discrete Poisson hardcore 1D model and  reinfections]
{\footnotesize Discrete Poisson hardcore 1D model and  reinfections}

\author[Ivan Cherednik]{Ivan Cherednik $^\dag$\\
}

\comment{
\authorcontributions{I.Ch.: research, data, exposition.}
\authordeclaration{no competing interest. }
\equalauthors{\textsuperscript{1} I.Ch. 
is the author this work.}
\correspondingauthor{\textsuperscript{2}
chered@email.unc.edu}

}

\begin{document}

\renewcommand{\tilde}{\widetilde}
\renewcommand{\hat}{\widehat}

\newcommand{\BR}{{\mathbb R}}
\newcommand{\BQ}{{\mathbb Q}}
\newcommand{\BC}{{\mathbb C}}
\newcommand{\BP}{{\mathbb P}}
\newcommand{\BZ}{{\mathbb Z}}
\newcommand{\BN}{{\mathbb N}}
\newcommand{\BS}{{\mathbb S}}

\newcommand{\cH}{{\mathcal H}}
\newcommand{\cA}{{\mathcal A}}
\newcommand{\cB}{{\mathcal B}}
\newcommand{\ccF}{{\mathfrak F}}
\newcommand{\cD}{{\mathcal D}}
\newcommand{\cL}{{\mathcal L}}
\newcommand{\cF}{{\mathcal F}}
\newcommand{\cP}{{\mathcal P}}
\newcommand{\cX}{{\mathcal X}}
\newcommand{\cY}{{\mathcal Y}}
\newcommand{\cS}{{\mathcal S}}
\newcommand{\cSol}{\hbox{$\mathcal Sol$}}
\newcommand{\cT}{\hbox{$\mathcal T$}}

\newcommand{\Z}{{\mathbb Z}}
\newcommand{\Q}{{\mathbb Q}}
\newcommand{\N}{{\mathbb N}}
\newcommand{\C}{{\mathbb C}}
\newcommand{\R}{{\mathbb R}}
\newcommand{\X}{{\mathbb X}}
\newcommand{\Y}{{\mathbb Y}}

\newcommand{\CH}{{\mathcal H}}
\newcommand{\CA}{{\mathcal A}}

\def\HH{\mbox{${\mathcal H}$\kern-5.2pt${\mathcal H}$}}

\newcommand{\binomial}[2]{\genfrac{(}{)}{0pt}{}{ #1 }{ #2 }}
\newcommand{\qbinomial}[2]{\genfrac{[}{]}{0pt}{}{ #1 }{ #2 }_q }
\newcommand{\qbinom}[3]{\genfrac{[}{]}{0pt}{}{ #1 }{ #2 }_{ #3 } }


\def\der{\partial}
\def\tensor{\otimes}
\def\gam{\gamma} \def\Gam{\Gamma}
\def\del{\delta} \def\Del{\Delta}
\def\kap{\kappa}
\def\lam{\lambda} \def\Lam{\Lambda}
\def\Comp{{\mathbb C}}
\def\sM{{\mathcal M}}

\newtheorem{theorem}{Theorem}[section]
\newtheorem{maintheorem}[theorem]{Main Theorem}
\newtheorem{proposition}[theorem]{Proposition}
\newtheorem{definition}[theorem]{Definition}
\newtheorem{lemma}[theorem]{Lemma}
\newtheorem{corollary}[theorem]{Corollary}
\newtheorem{notation}[theorem]{Notation}
\newtheorem{remark}[theorem]{Remark}
\newtheorem{example}[theorem]{Example}

\newtheorem{theorem }{Theorem}[section]
\newtheorem{maintheorem }[theorem]{Main Theorem}
\newtheorem{proposition }[theorem]{Proposition}
\newtheorem{definition }[theorem]{Definition}
\newtheorem{lemma }[theorem]{Lemma}
\newtheorem{corollary }[theorem]{Corollary}
\newtheorem{notation }[theorem]{Notation}
\newtheorem{remark }[theorem]{Remark}
\newtheorem{example }[theorem]{Example}

\newtheorem{ maintheorem }[theorem]{Main Theorem}
\newtheorem{ theorem}{Theorem}[section]
\newtheorem{ proposition}[theorem]{Proposition}
\newtheorem{ definition}[theorem]{Definition}
\newtheorem{ lemma}[theorem]{Lemma}
\newtheorem{ corollary}[theorem]{Corollary}
\newtheorem{ notation}[theorem]{Notation}
\newtheorem{ remark}[theorem]{Remark}
\newtheorem{ example}[theorem]{Example}

\newtheorem{thm}{Theorem}[section]
\newtheorem{prop}[thm]{Proposition}
\newtheorem{lem}[thm]{Lemma}
\newtheorem{cor}[thm]{Corollary}
\newtheorem{conj}[thm]{Conjecture}
\newtheorem{con}[thm]{Conjecture}
\newtheorem{dfn}[thm]{Definition}
\newtheorem{df}[thm]{Definition}
 \newcommand{\rem}{{\bf Comment.\ }}
 \newcommand{\rmk}{{\bf Comment.\ }}
 \newcommand{\exmp}{{\bf Example.\ }}
 \newcommand{\ex}{{\bf Example.\ }}
 \newcommand{\prob}{{\bf Problem.\ }}

\newtheorem{note}{Note} 
\renewcommand{\thenote}{}
\newtheorem{ack}{Acknowledgments}
\renewcommand{\theack}{}
\renewcommand{\appendixname}{\bf Appendix}

\hyphenation{
ap-pen-dix as-ymp-tot-ic at-trib-uted at-trib-ut-able
Bry-li-n-sky com-mu-ta-tion de-ge-ne-rate
de-riv-a-tive dis-trib-ute equi-vari-ant ex-tra-or-di-nary  
geo-met-ric griev-ance griev-ous grad-ed ho-lo-no-my ho-mo-thetic
in-fin-ite-ly in-fin-i-tes-i-mal Ha-rish Cha-n-dra mul-ti-plic-able 
non-euclid-ean non-iso-mor-phic non-smooth par-a-digm 
par-a-bol-ic pa-rab-o-loid pa-ram-e-trize phe-nom-e-non 
post-script pseu-do-dif-fer-en-tial pseu-do-fi-nite 
qua-drat-ics quad-ra-ture Han-kel rec-tan-gle semi-def-i-nite 
set-up wide-spread Euler-ian Feb-ru-ary Gauss-ian Grothen-dieck 
Hamil-ton-ian Her-mi-t-ian her-mi-t-ian Jan-u-ary 
Japan-ese Ka-shi-wa-ra Kor-te-weg Le-gendre No-vem-ber Rie-mann-ian 
Sep-tem-ber Za-mo-lo-d-chi-kov Kni-zh-nik quan-tum Op-dam
Mac-do-nald Ca-lo-ge-ro Su-ther-land Mo-ser 
Ol-sha-net-sky  Pe-re-lo-mov in-de-pen-dent ope-ra-tors 
cy-clo-to-mic ra-tio-nal de-gen-er-a-tion 
in-ter-est-ing de-for-ma-tions de-for-ma-tion pro-ce-dure 
fol-lows ope-ra-tors  pre-serve suf-fices ap-proach 
for-mu-las con-sider its com-ple-tion cor-re-spond-ing 
au-to-mor-phism be-cause pro-por-tional fi-nal-ly let-ting 
equi-v-a-lence ge-n-er-al-ized Mac-do-nald iden-ti-ties 
cor-re-s-pond sub-dia-grams par-ti-tion na-t-u-ral-ly 
or-dered stan-dard de-for-ma-tion ar-gu-ment com-bined 
sphe-r-i-cal rep-re-sen-ta-tions tri-go-no-me-t-ric
ge-n-er-al-ly speak-ing pri-m-it-ive ir-re-du-cible 
sum-ma-tion  rep-re-sen-ta-tives pro-por-ti-o-na-li-ty
ultra-sphe-ri-cal Ro-gers}

\def\ffor{\quad\hbox{ for }\quad}
\def\wwhen{\quad\hbox{ when }\quad}
\def\wwhere{\quad\hbox{ where }\quad}
\def\aand{\quad\hbox{ and }\quad}
\def\for{\  \hbox{ for } \ }
\def\iif{ \ \hbox{ if } \ }
\def\when{ \ \hbox{ when } \ }
\def\where{\  \hbox{ where } \ }
\def\and{\  \hbox{ and } \ }
\def\and{\  \hbox{ and } \ }
\def\oor{\  \hbox{ or } \ }
\def\proof{{\em Proof. \  }}

\def\equal{\stackrel{\,\mathbf{def}}{= \kern-3pt =}}

\def\la{\lambda}
\def\La{\Lambda}
\def\om{\omega}
\def\Om{\Omega}
\def\Th{\Theta}
\def\th{\theta}
\def\al{\alpha}
\def\be{\beta}
\def\ga{\gamma}
\def\ep{\epsilon}
\def\up{\upsilon}
\def\Up{\Upsilon}
\def\de{\delta}
\def\De{\Delta}
\def\ka{\kappa}
\def\kapp{\hbox{\bf \ae}}
\def\si{\sigma}
\def\Si{\Sigma}
\def\Ga{\Gamma}
\def\ze{\zeta}
\def\io{\iota}
\def\bio{b^\iota}
\def\aio{a^\iota}
\def\twio{\tilde{w}^\iota}
\def\hwio{\hat{w}^\iota}
\def\gio{\g^\iota}
\def\Bio{B^\iota}

\def\del{\delta}
\def\pa{\partial}
\def\vp{\varphi}
\def\ve{\varepsilon}
\def\inf{\infty}

\def\vph{\varphi}
\def\vps{\varpsi}
\def\vPh{\varPhi}
\def\vep{\varepsilon}
\def\vpi{{\varpi}}
\def\vth{{\vartheta}}
\def\vsi{{\varsigma}}
\def\vrh{{\varrho}}

\def\bph{\bar{\phi}}
\def\bsi{\bar{\si}}
\def\bvp{\bar{\varphi}}

\newcommand{\bS}{{\mathbf S}}
\newcommand{\bH}{{\mathbf H}}
\newcommand{\bF}{{\mathbf F}}
\newcommand{\bE}{{\mathbf E}}

\def\tal{\tilde{\alpha}}
\def\tbe{\tilde{\beta}}
\def\tde{\tilde{\delta}}
\def\tpi{\tilde{\pi}}
\def\txi{\tilde{\xi}}
\def\tPi{\tilde{\Pi}}
\def\tPhi{\tilde{\Phi}}
\def\tV{\tilde{V}}
\def\tJ{\tilde{J}}
\def\tla{\tilde{\lambda}}
\def\tga{\tilde{\gamma}}
\def\tGa{\tilde{\Gamma}}
\def\tvs{\tilde{{\varsigma}}}
\def\tu{\tilde{u}}
\def\tU{\tilde{U}}
\def\tw{\widetilde w}
\def\tW{\widetilde W}
\def\tB{\tilde B}
\def\tv{\tilde v}
\def\tV{\tilde V}
\def\tz{\tilde z}
\def\tb{\tilde b}
\def\ta{\tilde a}
\def\tih{\tilde h}
\def\trh{\tilde {\rho}}
\def\tx{\tilde x}
\def\tf{\tilde f}
\def\tg{\tilde g}
\def\tG{\tilde G}
\def\tk{\tilde k}
\def\tl{\tilde l}
\def\tL{\tilde L}
\def\tD{\tilde D}
\def\tR{\tilde R}
\def\tP{\tilde P}
\def\tH{\tilde H}
\def\tp{\tilde p}

\def\hH{\hat{H}}
\def\hh{\hat{h}}
\def\hR{\hat{R}}
\def\hY{\hat{Y}}
\def\hX{\hat{X}}
\def\hP{\hat{P}}
\def\hT{\hat{T}}
\def\hV{\hat{V}}
\def\hG{\hat{G}}
\def\hF{\hat{F}}
\def\hw{\widehat{w}}
\def\hW{\widehat{W}}
\def\hu{\hat{u}}
\def\hs{\hat{s}}
\def\hv{\hat{v}}
\def\hb{\hat{b}}
\def\hB{\widehat{B}}
\def\hze{\hat{\zeta}}
\def\hsi{\hat{\sigma}}
\def\hrh{\hat{\rho}}
\def\hth{\hat{\theta}}
\def\hy{\hat{y}}
\def\hx{\hat{x}}
\def\hz{\hat{z}}
\def\hg{\hat{g}}
\def\he{\hat{e}}
\def\hE{\widehat{E}}

\def\B{\mathbf{B}}
\def\I{\mathbf{I}}
\def\P{\mathbf{P}}
\def\G{\mathbf{G}}
\def\S{\mathbf{S}}
\def\F{\mathbf{F}}
\def\one{\mathbf{1}}
\def\Sn{\mathbf{S}_n}
\def\0{\mathbf{0}}
\def\H{\mathbf{H}}
\def\V{\mathbf{V}}

\def\f{\mathcal{F}}
\def\çF{\mathcal{F}}
\def\o{\mathcal{O}}
\def\t{\mathcal{T}}
\def\r{\mathcal{R}}
\def\l{\mathcal{L}}
\def\m{\mathcal{M}}
\def\k{\mathcal{K}}
\def\n{\mathcal{N}}
\def\d{\mathcal{D}}
\def\p{\mathcal{P}}
\def\cP{\mathcal{P}}
\def\a{\mathcal{A}}
\def\h{\mathcal{H}}
\def\c{\mathcal{C}}
\def\y{\mathcal{Y}}
\def\e{\mathcal{E}}
\def\v{\mathcal{V}}
\def\z{\mathcal{Z}}
\def\x{\mathcal{X}}
\def\s{\mathcal{S}}
\def\g{\mathcal{G}}
\def\u{\mathcal{U}}
\def\w{\mathcal{W}}
\def\i{\mathcal{I}}
\def\j{\mathcal{J}}
\def\b{\mathcal{B}}

\def\lan{\langle}
\def\llb{(\!(}
\def\ran{\rangle}
\def\rrb{)\!)}
 \def\dim{{\hbox{\rm dim}}_{\mathbb C}\,}
\def\lng{\hbox{\rm{\tiny lng}}}
\def\sht{\hbox{\rm{\tiny sht}}}
\def\sph{\hbox{\rm{\tiny sph}}}
\def\inv{\hbox{\rm{\tiny inv}}}

\def\br#1{\langle #1 \rangle}

\def\rank{\hbox{rank}}
\def\gl{\mathfrak{gl}_N}

\newcommand{\Aut}{\operatorname{Aut}}
\newcommand{\Hom}{\operatorname{Hom}}
\newcommand{\End}{\operatorname{End}}
\newcommand{\Ind}{\operatorname{Ind}}
\newcommand{\ad}{\operatorname{ad}}
\newcommand{\pr}{\operatorname{pr}}
\newcommand{\aweyl}{\tilde{\mathbb S}_n}
\newcommand{\hec}{{\mathcal H}^t_n}
\newcommand{\Func}{{\mathcal F}({\mathbb C}^n,{\mathcal H}^t_n)}
\newcommand{\tr}{\operatorname{tr}}
\newcommand{\Out}{\operatorname{Out}}
\newcommand{\Rad}{\operatorname{Rad}}
\newcommand{\Spec}{\operatorname{Spec}}
\newcommand{\id}{\operatorname{id}}
\newcommand{\Int}{\operatorname{Int}}
\newcommand{\ct} {\operatorname{ct}}

\newcommand{\rat}{{\mathbb Q}}
\newcommand{\real}{{\mathbb R}}
\newcommand{\cplx}{{\mathbb C}}
\newcommand{\zint}{{\mathbb Z}}

\newcommand{\sq}{\phantom{1}\hfill$\qed$}
\newcommand{\Rea}{\Re}
\newcommand{\Ima}{\Im}

\newcommand{\st}{\bowtie}
\newcommand{\modd}{\mbox{\,mod\,}}
\newcommand{\lr}{\langle}
\newcommand{\rr}{\rangle}
\newcommand{\eps}{\varepsilon}
\newcommand{\phk}{\phi^{(k)}}
\newcommand{\psk}{\psi^{(k)}}
\newcommand{\Res}{\mbox{Res}\;}
\newcommand{\sgn}{\mbox{sgn}}
\newcommand{\mn} {\left\{ \begin{array}{c}m\\
n\end{array}\right\}}

\def\sX{\mathscr{X}}
\def\sH{\mathscr{H}}
\def\sY{\mathscr{Y}}
\def\TT{\mathfrak{T}}
\def\JJ{\mathfrak{J}}
\def\HH{\mathfrak{H}}
\def\FF{\mathfrak{F}}
\def\GG{\mathfrak{G}}
\def\CC{\mathfrak{C}}
\def\LL{\mathfrak{L}}

\def\BB{\mathfrak{B}}
\def\AA{\mathfrak{A}}
\def\ZZ{\mathfrak{Z}}
\def\HH{\hbox{${\mathcal H}$\kern-5.2pt${\mathcal H}$}}
\def\HHH{\hbox{${\mathbb H}$\kern-4.2pt${\mathbb H}$}}
\def\tHH{\widetilde{\HH\ }}

\font\smm=msbm10 at 12pt 
\def\symbol#1{\hbox{\smm #1}}
\def\lsmash{{\symbol n}}
\def\rsmash{{\symbol o}}
\def\#{\sharp}

\font\tenbf=cmbx10
\font\tenrm=cmr10
\font\tenit=cmti10
\font\ninebf=cmbx9
\font\ninerm=cmr9
\font\nineit=cmti9
\font\eightbf=cmbx8
\font\eightrm=cmr8
\font\eightit=cmti8
\font\sevenrm=cmr7
\font\sevenbf=cmbx7


\begin{abstract}
We suggest a new hardcore Poisson-type 
distribution for Young diagrams with
the row lengths from some finite list. 
A discrete variant of the 
time-ordered Mat\'{e}rn II 
process in 1D is employed. This approach is related to that
based on the interlacing sequences due to Kerov
and others, but we restrict the number of rows.
The basic lengths
are assumed comparable with the total order of the diagram
in the quasi-classical limit, which results in new methods
and new formulas. An interesting application
is to random walks where the steps are at the points satisfying
the classical Poisson distribution or our truncated
one. In the simplest case, one obtains the distribution 
$\{e^{-\al}I_r(\al), -\infty < r < \infty\}$ for the Bessel
$I$-functions, which provides some probabilistic interpretation
of its many properties. 
An immediate application of our truncated Poisson
distributions is to 
modeling reinfections in epidemics, which is of
obvious importance for the Covid-19 pandemic.
\end{abstract}

\thanks{$\dag$ \today. \ \ \ Partially supported by NSF grant
DMS--1901796.} 

\address[I. Cherednik]{Department of Mathematics, UNC
Chapel Hill, North Carolina 27599, USA\\
chered@email.unc.edu}

\maketitle

\smallskip

{\small
\centerline{{\bf MSC} (2010): 05A15, 05A18, 33C10, 60C05, 60G55, 
60E05, 62P10,91F99}
}

\section{\bf Introduction}
The main general result of this paper 
is the definition and calculation of
a new Poisson-type 
distribution for sequences of  non-overlapping
subsegments of lengths $\{L_1,\cdots,L_s\}$ in a segment of 
length $N$.
Equivalently, this distribution is
for Young diagrams of order $\le N$ with the 
the row lengths $\{L_i\}$. This is
related to the interlacing sequences due to Kerov
and others, a starting 
point for many far-reaching applications: Jack-Whittaker-Macdonald 
functions and the corresponding stochastic point processes. 
From this perspective, we 
assume that the distances between the points are $\ge L_i$.  

The creation of these segments is
governed by a discrete variant of the time-ordered 
Mat\'{e}rn II process in 1D.
We assume that
that $\lim_{N\to \infty} L_i/N=\nu_i>0$ in the quasi-classical
limit, so the 
edge effects are significant. 
Taking the limit required a procedure
which seems new even for one $L$, as well as  
the resulting formulas.

A variant is the probability
distribution for random walks with the steps $\pm 1$ at the points 
satisfying the Poisson $\al$-distribution 
or our truncated one. In the simplest case, we arrive at 
$\{e^{-\al}I_r(\al), -\infty\! <\! r\! <\! \infty\}$ for the 
Bessel $I$-functions, which provides a  probabilistic 
interpretation of the properties of $I$-functions.
The restriction of this approach 
to Catalan (non-negative) paths is considered.  

\vskip 0.2cm
An immediate application of our truncated Poisson
distribution is to 
modeling reinfections in epidemics, aimed at  
the Covid-19 pandemic.
Combinatorially, segments of one or several 
different lengths are protective immunity intervals; they
are placed in a bigger segment, the epidemic cycle. Practically, 
no greater than 2-3 different strains can be present simultaneously. 
The case of 1 strain (when $s=1$) is the key.

\vskip 0.2cm
Reinfections were relatively rare before Covid-19 (during 1 cycle). 
The Poisson distribution can be
expected if the immunity and the duration of the disease 
are disregarded. 
However, immunity is the key here. There were
not many papers on modeling reinfections; see e.g. \cite{ADDP},
which was SIR-based. We mention that the SIR-models proved to be 
not applicable to Covid-19; see \cite{ChB, ChM}.  
\vskip 0.2cm
\vfil

{\bf Hardcore point processes.}
The corresponding mathematical tool is the theory of
{\sf Poisson hardcore point processes}, more specifically,
the time-ordered {\sf Mat\'{e}rn process II}  and in its {\sf 
lattice} variant. There are
quite a few processes where the distances between neighboring
objects must be greater than some
constant. 
The usual examples are {\sf forestry, ecology, vehicular networks,
cellular networks}, etc. Also, see \cite{Gi} for {\sf Tonks gas}. 
In statistical physics,  ``small systems" are of this type, those
far from thermodynamic equilibrium. 
There is a vast literature on  
Mat\'{e}rn processes I,II,III, mostly in 2D. 
See e.g. \cite{KD} on vehicular networks and references there. 
Modeling {\sf vehicular networks}, clear 1D processes, is
somewhat similar to modeling reinfections. The  
interval between cars is a counterpart of the
immunity intervals. 

\vskip 0.2cm
In such and similar
examples, the intervals between objects
are mostly assumed small vs. the
domains where they are considered. 
Accordingly, the {\sf edge effects} are mostly ignored. 
This is different for epidemics: 
the immunity intervals  are quite comparable with the 
duration of the  epidemic cycle.  The continuous 
distribution we 
obtain in  (\ref{distprim}) is a certain
truncated version of the classical {\sf Poisson distribution};
we think it is new. It is not its straightforward  truncation via the
relative probabilities.


\vskip 0.2cm

We begin with the following discrete
setting: configurations of subsegments of the length $L\!+\!1$ 
with gaps between them in a given segment
of the length $N$. This can be necessary 
when the process depends on $N$ and $L$,
not just on $L/N$ in the limit (see below).
In the lattice version,
exact formulas can be obtained in terms of binomial coefficients
and generalized to subsegments of any lengths, which 
requires only basic combinatorics. We make them as explicit
as possible in Corollary \ref{cor:pir}. 
\vfil

{\sf Edge effects.} This is when the last segment
goes beyond $N$. They naturally result in a sum
of $L+1$ binomial coefficients. 
To perform the continuous
limit, we need a 
formula where the number of terms depends only
on the number of subsegments.  Generally, the assumption is that 
the size of the subsegments
is comparable with that of
the segment where they are considered. 
For one $L$,  the  continuous limit is when 
$\lim _{N\to \infty} L/N=\nu>0$.  The resulting
distribution has  Ceiling$[1/\nu]$ states. The process
of finding  this limit suggested in the paper is of interest;
the procedure we use seems new.   
\vfil

\vskip 0.2cm
{\sf Several lengths.}
A natural extension is from one $L$ to $\{L_i, 1\le i\le s\}$
(with different probabilities). We do this
in 2 stages: first, for $s=2$ and then the general case.
Interestingly, the formulas remain quite compact, and the
limiting procedure for one $L$ can be used almost without 
modifications.
See Theorems \ref{thm:Fdif-new} and \ref{thm:lim-new}.
The probability distribution becomes that for the Young diagrams 
with $\le N$ boxes
and the rows of sizes from  $\{L_i\}$. Thus, we restrict
the number of rows in contrast to \cite{Ke,BO,Ol} and other papers.
Also, the creation of the segments (rows in the 
Young diagrams) is subject to the Poisson-type distribution.  
Generally, one can explicitly calculate the correlation functions
for these processes; we provide the simplest ones.
\vskip 0.2cm
\vfil

{\sf Bessel-type formulas.}
A variant of our approach is when we consider $L_i$ as some
jumps of the energy function or similar functions
with some probabilities.
It is assumed that the energy constantly increases by $1$ unless 
for the jumps. It is natural here to allow the jumps by $\pm L_i$,
not only by $+L_i$ as for reinfections and similar processes.
In the quasi-classical limit $\lim_{N\to \infty} L_i/N=\nu_i$,
it leads to interesting multi-dimensional Bessel-type functions.

For one $L$ and when $\nu\to 0$, i.e. when the segments become
points, we obtain a random walk with
probabilities $p,q$ of the steps $\pm1$ that occur at the
points subject
to the Poisson $\al$-distribution. The probability of
$r$ outputs is then $e^{-\al}I_r(\al)$ for $p\!=\!\frac{1}{2}\!=\!q$.
Interestingly, {\sf pure} Bessel functions $I_r$ occur here. 
It is reasonable here to allow only positive $r$ during the process,
i.e. to impose the Catalan condition for the paths; we provide the
formulas.

\vskip 0.2cm

{\bf Reinfections (Covid-19).}   
The duration of the  Covid-19 epidemic (from late 2019)
is already beyond 2 years; we are still in its 1st cycle.
This epidemic  
was practically uninterrupted except for minor breaks between
the waves (mostly during summer periods). Due to 
the unusually large  number of the strains of Covid-19, all with 
with very high transmissibility, the natural immunity did not last
too long for Covid-19, as well as the
immunity due to the vaccinations. 

For instance, those infected
by the ``wild strain" 
(the G-strain dominated in Europe in early 2020) could be reinfected
by Alpha, then by Delta  
(B.1.617.2 and AY lineages), and then by Omicron 
(B.1.1.529 and BA lineages). The average immunity durations were 
not that long, presumably
about 5-8 months; the waves of different strains obviously 
contributed to this. We will disregard below relatively rare cases
when someone is infected simultaneously by 2 strains.


\vskip 0.2cm
The statistics of reinfections and the recurrences for Covid-19 
is not very reliable. Some countries reported only the
total number of (known) infected individuals, not the total
number of detected infections.
For instance, this was the case with England until 
January 31, 2022; the data on Covid-19 in England
are generally among the most systematic.  According to the
UK Health Security Agency (UKHSA), the 
number of detected reinfections can be about
10\% in early 2022. 
In quite a few countries, it was significantly greater than this and
reinfections were present well before 2022. 
\vfil

We note that the available data are for the detected cases.
Massive testing began at the end of 2021 
in quite a few countries, but very many cases remain unreported.
However even the detected cases of
{\sf double} reinfections (3 Covid-19 infections)
appeared not too rare. The actual numbers can be significantly
higher among all infections, including asymptomatic ones. 
\vfil

We use in this paper a general method, which
can be naturally extended to any number of {\sf parallel 
infections} with different immunity durations of any
any lengths. There is 
only one constraint: 
the corresponding subsegments must not overlap.

We note that  Covid-19 provides quite a few examples of parallel strains:
\{Alpha and Delta\}, \{Delta: AY.1, \ldots, AY.4.2\},
\{Delta, Omicron\}, 
\{Omicron: BA.1, BA.2, BA.3, \ldots\}. It is really 
rare when someone could be infected by two parallel strains 
at the same time. 

\vskip 0.2cm

{\sf Main hypotheses.}
We assume that people are exposed 
to the infection {\sf uniformly} during the cycle of the epidemic,
which is $N$ days in the paper, with 
probability $\be$ per day. According to \cite{ChB,ChM},
the curves of the total number of detected infections
in very many countries (all we considered)
are essentially of Bessel type for phase 1 and are of 
linear type for phase 2; a clear phase transition can be
seen in many countries. However, it is very reasonable
to assume that statistically the process
is not far from linear. 

Moreover, the numbers of consecutive waves in many countries
was like 3-6, which provides another reason
to assume that the spread of Covid-19 is
linear {\sf statistically}.
Generally, the {\sf Law of Large Numbers} (LLN) is always a rationale
for the {\sf uniformity assumption}. The averages over
3-6 waves are sufficient for this. Similarly,  
not much will change statistically for  reinfections 
if we try to incorporate the ``exact" (Bessel-type) shapes 
of the curves of the total cases.   
\vskip 0.2cm
\vfil

The 2nd hypothesis is that the 
impact of the vaccinations and (significant) 
number of undetected and  asymptomatic cases can be addressed
via diminishing  the {\sf susceptible population}
of a country. The size of population does
not directly appear in the formulas. The vaccinations
generally decrease $\be$, but the strains of
Covid-19 increased their transmissibility  
during the epidemics. 

The 3rd hypothesis is a minor one: we disregard the
duration of the disease. It is simply added to the
immunity interval. 

We think that all 3 assumptions are
quite reasonable for the last 2 years of Covid-19.
So the challenge is to provide the distribution of
the reinfections based on them, and then adjust it
to the real data. 
\vskip 0.2cm
\vfil

{\sf The distribution}. Under
these assumptions, we define and calculate the probabilities 
$\pi_r$ for
$r$ infections during the total period.
They depend on 
the duration of an epidemic, 
$N$ days, the immunity interval $L$, and the probability $\be$ to
be exposed to the infection during 1 day. Later, we make
$\be=\frac{\al}{N}$ for some parameter $\al$. We note that it can be
only among reported case; this does not influence our
analysis but $\al$ and $L/N$ can change.

The parameters $\al,\be$ can be  determined if the total number
of non-infected people (during the whole cycle) and those
with exactly 1 infection are known.
In the absence of immunity,
$\al$  can be estimated as the total number of infected people
divided by the size of the population in the area.
Generally, $e^{-\al}$ is the 
number of {\sf noninfected people} (during the whole cycle
of epidemic) divided by the 
population of the country. This is 
as for the classical Poisson distribution. 

Then 
we (approximately) find $L$ using our formula
for $\pi_1$ or its limit $\pi'_1$
(exactly one infection during $N$ days). For the classical
Poisson distribution (no immunity), it is $\al e^{-\al}$.
Now it depends on $L$.

So we can assume that $\be=\al/N$ and $L$ are known.
Then  $\pi_2$ and $\pi_3$ are the probabilities of 
$2$ and $3$ infections during
$N$ days, which can be compared with reinfection data.

We consider through the
paper 3 {\sf model situation}: $N=750$ (about 25 months),
$L=150$ (about 5 months), and 3 values of $\al=1, \log(2), 0.5$
in (\ref{p-val2}). 
For these values of $\al$, 
about $37\%, 50\%, 60\%$ of the susceptible 
population remain noninfected during 25 months. The 3rd case
basically matches the number of reinfections
for Covid-19 reported in England (until 02/2022). 
\vskip 0.2cm
\vfil

{\bf Main findings.}
The famous Poisson distribution is a straightforward
limit of a simple distribution in terms of binomial coefficients,
its lattice variant. 
Our distribution $\{\pi_r\}$ for
the probabilities of $r$ infections is a sum of $L\!+\!1$  
binomial coefficients, where $L$ is the immunity duration.
We calculate $\pi_r$ and $\pi'_r=\lim_{N\to \infty}\pi_r$ 
as $\lim_{N\to \infty} L/N=\nu>0$;
the limiting procedure is interesting.
We give other formulas for $\pi_r$ 
and those for the corresponding generating functions. 
\vskip 0.2cm

In spite of the 1D setup, there are various applications
of $\pi_r, \pi_r'$, not only for reinfections. 
Almost any networks have refractoriness: excited agents cannot
be immediately re-excited.  Vehicular networks and
trading equities in stock markets are typical examples.
We focus in this paper on networks with relatively small number
of possible states and when $\nu=L/N$ cannot be assumed negligible,
which is obviously the case with reinfections. 
\vskip 0.2cm

\comment{ 
Our Poisson-type distributions
can be applicable to other 1D Mat\'ern-II processes
when the lengths of (disjoint) random subsegments
are comparable with the total length of the considered segment
and when the last subsegment can go beyond $N$. 
The subsegments
can be of different types (with different lengths and probabilities), 
which links our paper to stochastic processes. General formulas
for any number of precesses and
the corresponding distributions of probabilities
of Young diagrams are obtained in (\ref{pipims}) and (\ref{piprims}). 
}

As far as we know, our  distributions
$\{\pi_r\}$ and $\{\pi'_r\}$ are new,  as well as their
application to reinfections of epidemics. 
The number of sates is bounded in our approach:\,
$r\le $Ceiling$[N/L]$. The edge effects are important;
we allow one of the subsegments  to go beyond $N$, the right
endpoint. This assumption
is necessary for reinfections.
Stock markets are such too: there can be
open positions after the end of the considered period.
Mathematically, the corresponding sum of all
probabilities will not be $1$ without the edge effects. 

This approach can be smoothly extended to any number
of lengths: $\{L_i,1\le i\le s\}$ with the
corresponding probabilities.
The formulas remain reasonably compact; the
limiting procedure is basically the same as for one $L$.
See Theorems \ref{thm:Fdif-new} and \ref{thm:lim-new}.
The probability distribution becomes that for the Young diagrams
of size $\le N$  
with finitely many rows: those of lengths $L_i$.

An interesting variant here
is when the ``energy function" can jump by $\pm L_i$ with
probabilities $p_i,q_i$ and increases by $1$ otherwise. 
This is motivated by physics, networks and share-prices. 
When $N\to \infty$ and $\nu_i\to 0$,
i.e. the segments $L_i$ become points, we arrive at the random
1D walk with jumps up and down at the points satisfying some
Poisson-type distribution or our truncated one. We arrive at
generalized Bessel $I$-functions. Catalan-type (non-negative) 
paths are of obvious interest here.
We provide only  Theorem \ref{prop:catal}
for $s=1$ (for one $L$).

\section{\bf Hardcore Poisson-type processes}
If the immunity factor is omitted, the distribution
of reinfections is as follows. Assume that an
epidemic lasts $N$ days and $\be=\al/N$ is the probability
to be infected during one day. Then the probability to
be infected $r$ times during $N$ days and its continuous limit
are given by the classical Poisson distribution and its 
combinatorial counterpart. Namely:
\begin{align}\label{P1}
p_r\!=\!\binom{N}{r}\be^r (1\!-\!\be)^{N-r} \text{\, for\, } 
r\!=\!0,1,2,\ldots\,,
\, p'_r\!=\!\lim_{N\to \infty} p_r\!=\!\frac{\al^r}{r!e^{-\al}},
\end{align}
where $e$ is the Euler number. 

{\sf Three basic examples.} Let $\al=1, \al=\log(2)\approx 0.69, 
\al=0.5$.  Then, 
{\small
\begin{align}\label{p-val}
 p'_0=&p'_1=\frac{1}{e}\approx 0.37, 
\, p'_2=\frac{1}{2e}\approx 0.18,
\, p'_3=\frac{1}{6e}\approx\, 0.06 \text{\, for\, }\al=1,\\
 p'_0=&0.5, \ p'_1=\log(2)/2 \approx 0.35,\ p'_2\approx 0.12,
\ p'_3\approx 0.03 \text{\, for\, } \al=\log(2),\notag\\
&\text{and \ } p'_0\approx 0.61, \ p'_1\approx 0.30,\ p'_2\approx 0.08,
\ p'_3\approx 0.01 \text{\, for\, } \al=0.5. \notag
\end{align} 
}
\!\!The corresponding values for the combinatorial
 $p_i$ are about the same.
\vskip 0.2cm

{\sf Adding immunity.}
Assume that an individual infected at day $x$ cannot be
infected again for days  $x+1,x+2,\ldots, x+L$, i.e.
$L<N$ is the duration of the immunity interval. 

Let $\pi_r$ be the probability
of $r$ infections during $N$ days for $r=0,1,\ldots\, .$
If $1\le x_1<x_2<x_3<\ldots<x_r\le N$ are the
infection days, then 
 $x_1<x_2-L<x_3-2L\ldots x_r-(r-1)L$ and there are 2 cases:
\vskip 0.2cm

\centerline{{\bf (a)}\  $x_r+L\le N$,\ \  and, otherwise, 
{\bf (b)}\  $x_r+L> N$. }
\vskip 0.2cm

Here $x_i$  are for the actual infections, when the
disease begins. The potential 
infections are  the days when an individual 
was exposed to the infection, which is assumed with probability $\be$.
Due to the immunity, not all of the exposures result in the
actual infection (disease).
Any number of potential infections can occur (anywhere)
 during the periods 
$x_i+1,\ldots, x_i+L$ for $1\le i< r$ and during 
the end  period $x_r+1,\ldots, \min\{x_r+L,N\}$. This means
that  these periods can be removed from the
consideration when counting the
probabilities. 
Switching to $x_i'=x_i(i-1)$ for $1\le i \le r$,
the probability of the event 
``$x'_1<x'_2<\ldots<x'_r$" is $\be^r (1-\be)^{N-rL-r}$
for {\sf (a)} and $\be^r (1-\be)^{(x_r-(r-1)L -r)}$ for {\sf (b)}.
\vskip 0.2cm

Let $P(r,N,\be)\equal \binom{N-Lr}{r}\be^r
(1-\be)^{N-Lr-r}$; it is $0$ if $r<0$ or when 
$N< Lr+r$.
We obtain the following straightforward formula:
\begin{align}\label{Pform}
&\pi_r=P(r,N,\be)+\be\sum_{i=1}^L P(r-1,N-i,\be), \ r=0,1,\ldots\ .
\end{align} 
To give an example: for $L=1$, 
$\pi_r=\be^r(1-\be)^{N-2r}\left(\binom{N-r+1}{r}+
\be\binom{N-r}{r}\right)$.
Here and for any fixed $L$, $\lim_{N\to \infty} \pi_r=p'_r$,
where $p'_r$ are from (\ref{P1}), where we set $\be=\al/N$. 

One has: $\sum_{r=0}^\infty \pi_r=1$, where $r\le \frac{N+L}{L+1}$
are sufficient in this sum. This is some combinatorial identity,
which immediately follows from the definition of $\pi_r$. 
Obviously,  $\pi_0=p_0$ for any $N,L$. 

For our three basic examples above, we will take
$N=750, L=150$. Then, $\nu=0.2$ and 
\begin{align}\label{p-val2}
&\pi_0\approx 0.37,\, \pi_1\approx 0.44, 
\, \pi_2\approx 0.17,
\, \pi_3\approx 0.02 &\text{\, for\, }\al=1,\\
&\pi_0\approx 0.50,\, \pi_1\approx 0.39,\, \pi_2\approx 0.10,
\, \pi_3\approx 0.01 &\text{\, for\, } \al=\log(2),\notag\\
&\pi_0\approx 0.61,\, \pi_1\approx 0.33, 
\, \pi_2\approx 0.06,
\, \pi_3\approx 0.004 &\text{\, for\, }\al=0.5.\notag
\end{align} 
The change is not dramatic vs.
(\ref{p-val}) since $\al$ and $\nu$ are relatively small.
For instance,  
$\pi_1\approx 0.51$ if $\nu=0.4$ for $\al=1$
(with the same  $\pi_0\approx 0.37$). 
\vskip 0.2cm

Let us rewrite the formula
for $\pi_r$ without the $L$-summation.

\begin{theorem}\label{thm:Fdif}
Let 
$F_r(X)=d^{r-1}\left(X^{N-Lr}\,\frac{\,1-X^L}{1\, -\, X}
\right)/dX^{r-1}$ (the
$(r-1)${\small th} derivative), where $r>0,L>0$ and $N-Lr>0$. Then
{\small
\begin{align}\label{difform}
&\pi_r=\binom{N-Lr}{r}\,\be^r\, (1-\be)^{N-Lr-r}+
\frac{\be^r}{(r-1)!}\, F_r(X\mapsto 1-\be). \hfill \text{\sq}
\end{align} 
}
\end{theorem}
\vskip 0.2cm
The number of terms in this formula depends only on $r$ 
(not on $L$), which is the key when
considering its limit as $L,N\to \infty$.
\vskip 0.2cm

\section{\bf Continuous limit}
Let us provide the first 4 cases of (\ref{difform}),
where we use directly the theorem:

{\small
\begin{align}
\pi_0= (1-\be)^N,\, 
\pi_1&=\be(N\!-\!L)\,(1-\be)^{N-L-1}+(1-(1-\be)^L)\,
(1-\be)^{N-L},\notag\\
\pi_2=\be^2\binom{N\!-\!2L}{2}&\,(1\!-\!\be)^{N-2L-2}
+\be (N\!-\!2L)
\left(1\!-\!(1\!-\!\be)^L\right)(1\!-\!\be)^{N-2L-1}\notag\\
+&\Bigl(1-(1-\be)^L-\be L (1-\be)^{L-1}\Bigr)(1-\be)^{N-2L},\notag\\
\pi_3=\frac{\be^3}{6} \Biggl(&(N\!-\!3 L) (N\!-\!3 L\!-\!1) 
(N\!-\!3 L\!-\!2) x^{N\!-\!3 L\!-\!3}\,+ \notag\\
3\,\biggl(-2 L x^{L-1}\, \bigl(\,&\frac{x^{N\!-\!3 L}}{\be^2}
+\frac{(N\!-\!3 L) x^{N\!-\!3 L\!-\!1}}{\be}\,\bigr)-\frac{L(L\!-\!1) 
x^{N\!-\!2 L\!-\!2}}{\be}\,+\notag\\
\left(1\!-\!x^L\right)\Bigl(\frac{2 x^{N\!-\!3 L}}{\be^3}+
&\frac{2 (N\!-\!3 L) x^{N\!-\!3 L\!-\!1}}{\be^2}+
\frac{(N\!-\!3 L) 
(N\!-\!3 L\!-\!1) x^{N\!-\!3 L\!-\!2}}{\be}
\Bigr)\!\!\biggr)\!\!\Biggr).\notag
\end{align}
}

Here we set $x\equal 1-\be$. Also,
$N\ge L$, $N\ge 2L$ and $N\ge 3L$ correspondingly; 
generally, $r\le \frac{N}{L}$.
Recall that $\pi_r>0$
for $N\ge L(r\!-\!1)+r$, i.e. for $r\le \frac{N+L}{L+1}$.
 For instance, $\pi_1=1-(1-\be)^N$ for
{\sf any} $L$ such that $L>N$. 

Notice that $\pi_0$ does not
depend on $L$; this is obvious because the duration of the
immunity interval does not affect those non-infected. 
\vskip 0.2cm

Setting $\be\equal\al/N$ and $\lim_{N\to \infty}
N/L= \nu\ge 0$, the limits $\pi'_r=\lim_{N\to \infty}\pi_r$
for $r=0,1,2$ are as follows: 
\begin{align}
\pi_0'=e^{-\al},\,\ 
\pi_1'=e^{(\nu-1)\al}\bigl(1+\al(1-\nu)\bigr)& - e^{-\al},\notag\\
\pi'_2=e^{(2\nu-1)\al}\left(1+\al(1-2\nu)+ 
\frac{\al(1-2\nu)^2}{2}\right)-
e&^{(\nu-1)\al}\left(1+\al(1-\nu)\right).\notag
\end{align}

\comment{
{\small
\begin{align*}
\pi'_3=&e^{3\nu-1}\left(1+\al(1-3v)+\frac{\al^2(1-3\nu)^2}{2}+
\frac{\al^3(1-3\nu)^3}{6}\right)\\
&-e^{2\nu-1}\left(1+\al(1-2v)+\frac{\al^2(1-2\nu)^2}{2}\right).
\end{align*}
}
}

We assume here  
that $r\nu \le 1$. Generally, under this assumption:
{\small
\begin{align}\label{distprim}
&\pi'_r=e^{(r\nu-1)\al}\left(1+\al(1-r\nu)+
\frac{\al^2(1-r\nu)^2}{2}+\ldots+
\frac{\al^r(1-r\nu)^r}{r!}\right)\\
&-e^{((r\!-\!1)\nu-1)\al}\left(1+\al(1-(r\!-\!1)\nu)+
\ldots+\frac{\al^{r-1}(1-(r\!-\!1)\nu)^{r-1}}{(r-1)!}
\right).\notag
\end{align}
}

The last value of $r$ here is $r_\flat=$Floor$[\frac{1}{\nu}]$, 
where Floor$[x]$ is the integer part of $x$. The initial
inequality $r\le \frac{N+L}{L+1}$ gives that the last 
nonzero $\pi'_r$ is for 
$r_\#=$Ceiling$[N/L]$, which is
$r_\flat+1$ if $\frac{1}{\nu}$ is not an integer, and $\pi_\flat$
if $1/\nu\in \Z_+$.   
Indeed, $\frac{N}{L}+1>\frac{N+L}{L+1}>\frac{N}{L}.$
When $r_\#=r_\flat+1$, one has:
{\small
$$\pi'_{r_\flat+1}=1- e^{(r_\flat\nu-1)\al}\left(1+\al(1\!-\!r_\flat\nu)+
\frac{\al^2(1\!-\!r_\flat\nu)^2}{2}+\ldots+
\frac{\al^{r_\flat}(1\!-\!r_\flat\nu)^{r_\flat}}{(r_\flat)!}\right). 
$$
}

The positivity of $\pi'_r$ for $0\le r\le r_\flat+1$
can be readily seen  
from these formulas, though it of course follows
from the origin of $\pi'_r$. Here we calculated $\pi'_{r_\flat+1}$
directly from the definition. Alternatively, it can be obtained 
from the identity $\sum_{r=0}^{r_\#} \pi'_r=1$. 
This sum is obviously $1$ (the telescoping summation), which
holds {\it a priori\,}  
because  $\sum_{r=0}^{r_\#} \pi_r=1$, which is  due to the definition
of $\pi_r$.  We arrive at the following theorem.

\begin{theorem}\label{thm:Fdif-lim}
Assume that $N\to \infty$, $\be =\al/N$,
$\lim_{N\to \infty} L/N=\nu$ for some $\al>0$ and $0\le \nu \le 1$.
Then formula (\ref{distprim}) holds for any $0\le r \le r_\flat$,
as well as the additional formula above for $r_\#=r_\flat+1$ when 
$\frac{1}{\nu}$ is not an integer.
\end{theorem}
{\it Proof.} We can simplify 
$F_r$ from (\ref{difform}) 
considered in the limit as follows.
Let $\Phi_r\equal
\b^{r-1}\bigl(\frac{X^{N-Lr}-X^{N-L(r-1)}}{1-X}\bigr)$
for the  following 
{\sf formal\,} differentiation $\b$ of the ring generated by $X^{N-sL}$
and $1/(1-X)$ {\sf treated as independent symbols}\,:
{\small
$$
\b^p(X^{N-sL})\!=\!(1\!-\!s\nu)^p X^{N-sL}, \ 
\b^p\Bigl(\frac{1}{1-X}\Bigr)\!=\!
\frac{\be^p\, p!}{(1-X)^{p+1}},\  p=0,1, 2,\ldots \ .
$$
}
\!\!This differentiation is  $\be\, d/dX$ with the following
simplifications 
due to taking the limit.
First, we replace $\binom{N- s L}{r}$ 
by $\frac {N^r}{r!}$
for $s\le r$ due to $r\ll N-sL$. Second, 
 $X^{N-sL\,\pm p}$ is replaced  
by $X^{N-sL}$ for $0\le p\le r$ because these powers will
be finally evaluated at $X=1-\be=1-\frac{\al}{N}\to 1$.

We obtain that 
$\pi'_r=\frac{\al^r(1-r\nu)^r}{r!}\,e^{(r\nu-1)\al}
+ \frac{1}{(r-1)!}\,\Phi'_r$,
where $\Phi'_r$ is the limit of $\be\Phi_r$ after
the
evaluation $1\!-\!X\mapsto \be,\, X^{N-Ls}\mapsto e^{\al(s\nu-1)}.$

One has:\, $\Phi_r=\sum_{s=0}^{r-1}\binom{r-1}{s}\,
\b^s\Bigl(X^{N-Lr}-X^{N-L(r-1)}\Bigr)\, \b^{r-1-s}\Bigl(
\frac{1}{1-X}\Bigr)=$
{\small
\begin{align*}
&\sum_{s=0}^{r-1}\be^{r-s-1}\,\frac{(r-1)!}{s!}\,
\biggl(
\frac{\al^s (1-r\nu)^s X^{N-Lr}-\al^s (1-(r-1)\nu)^s X^{N-L(r-1)}}
{(1-X)^{r-s}}\biggr).
\end{align*}
}
Finally, $\frac{\be^{r-s-1}}{(1-X)^{r-s}}\mapsto
\be^{-1}$,  $\lim_{N\to \infty}(1\!-\!\be)^{N-Ls}= e^{\al(s\nu)-1}$,
and  
{\small
$$ \Phi'_r\!=\!
\sum_{s=0}^{r-1}\frac{(r\!-\!1)!}{s!}\,
\Bigl(\al^s (1\!-\!r\nu)^s e^{\al(r\nu-1)}- 
\al^s (1\!-\!(r\!-\!1)\nu)^s e^{\al((r-1)\nu-1)}\Bigr).
$$
}
\!\! This gives formula (\ref{distprim}). The calculation
is similar for $\pi_{r_\#}$. \sq
\vskip 0.2cm

The $\nu$-dependence of $\pi'_r$ is interesting.
For instance, since $\pi'_0$ does not depend on $\nu$,
$\pi'_1$ increases if $\nu$ increases. Indeed, the chances of
reinfections (counted by $\pi'_r$ for $r\ge 2$) diminish. 
Similarly, $\pi_{r_\#}$
decreases if present (if $1/\nu$ is not an integer).
Generally, we have the following
straightforward corollary.

\begin{corollary}\label{cor:der}
Let $\de_r=r\al\,e^{(r\nu-1)\al}\,\frac{\al^r(1-r\nu)^r}{r!}$
for $0\le r\le r_\flat$, and $\de_{-1}=\de_{r_\#}=0$. Then
$d\,\pi'_r/d \nu=\de_r-\de_{r-1}$ for $0\le r\le r_\#$.
In particular, 
$d\,\pi'_r/d \nu\ge 0$ for $1\le r\le \frac{1}{\nu}$
if and  only if\, $\nu\al\, e^{\nu\al}\ge (r-1)(1+z)^{r-1}z$ for
$z=\frac{\nu}{1-r\nu}$. Otherwise, this derivative is negative.
\sq
\end{corollary}

Practically, triple reinfections ($r\!=\!4$)
are hardly possible for one cycle of
any epidemic.  Though there can
be other random processes of this kind where
big $r$ make sense.  The distribution in 
(\ref{distprim}) is some quantization of
the Poisson distribution, where $\nu\to 0$ is the quasi-classical
limit. 
Accordingly, (\ref{Pform}) is its ``quantization" with $2$ parameters,
$L$ and $N$. 
One more parameter can be added to (\ref{Pform})
by switching to the $q$-binomial coefficients there, which we
will not discuss. 
\vskip 0.2cm

{\sf Reinfections in England.} As a demonstration,
let us try to employ these formula to the Covid-19
data from England.
``As of 31 January (2022), updated figures for England show 14845382 
episodes of infection since the start of the pandemic with 
588114 (4.0\%) reinfections covering the whole pandemic."
So, approximately $14845382-588114=14257268$ people were detected
to be infected at least once. Let us assume conditionally that
about $35$M were involved in collecting the data; the population of
England is about $57$M. Our approach can be applied
if only {\sf detected\,} cases and reinfections 
are taken into account; however, $\al$ depends on the number
of {\sf all\,} infections, including the asymptomatic and
undetected ones. Technically, we diminish $57$M to
$35$M, but this can be done directly via $\al$ (to adjust 
$\pi_0, \pi_1$). 
  
As in the 3 basic cases, we take $N=750$ and $L=150$. Then  
$\al\approx 0.5$; indeed, $\pi'_0=
e^{-0.5}\approx 0.6\approx 1-14/35$.  This is
basically the 3rd case in (\ref{p-val2}):  
$\pi_1\approx 0.33$, $\pi_2\approx 0.06$. 
Qualitatively, $0.06$ matches the data from UKHSA: about $0.04$
for $\pi_2$ (until January 31, 2022). 

\comment{
Physically, this is a uniform propagation of particles
of size $L$ on a segment of size $N$:
after a particle is created, its area cannot be used for
further particles. If $\nu=L/N$ is not too small,
then (\ref{Pform}) is supposed to be used. 
This can be with droplets  of liquid of fixed size
under some condensation process as particles, but then
the dimension is $2,3$.
}
\vskip 0.2cm

\section{\bf Poisson-Catalan distribution} Let us assume that a system 
is subject
to $N$ random events of $3$ sorts: (i)\, adding $1/N$ to its 
{\sf\em energy} $E$
with probability $(1-\frac{\al}{N})$, (ii)\, adding $L/N$ to $E$ with 
probability $\frac{\al p}{N}$ and (iii)\, subtracting $L/N$ with
probability $ \frac{\al q}{N}$, where  $p,q\ge 0$ and $p+q=1$.
This can be a ``slow" linear growth of energy with more significant 
accidental  transitions to excited states and back.  Another
interpretation is when $E$ is 
a share-price subject to some constant trend  with
relatively rare $\pm L/N$ fluctuations.

We begin with $E=0$ and fix the final energy
in the range $1\le E<1+L/N$. It is allowed for $E$ to take any
non-negative values before the final point but
the intermediate balances of additions 
and subtractions of $L/N$
must be always non-negative (the gains due to other points
are disregarded).  I.e. it must be a {\sf
Catalan path}.

We assume that $\lim_{N\to \infty} L/N=\nu\le 1$
and consider below only this limit. Accordingly,
$\p'_{r,m}$ be the probability that $E$
changes from $0$ to the final value in the range $0\le E<1+\nu$ for
Catalan paths 
with $r+m$ events of adding $\nu$ and $m$ 
events of subtracting $\nu$.

In the following proposition,  $\pi'_{r,m}$ denotes
 $\pi'_{r+2m}$ calculated 
for the parameters
 $\al_m=\text{\footnotesize $(1\!+\!2m\nu)$}\,\al$ 
and $\nu_m=\frac{\nu}{1+2m\nu}$.
In fact, only $\al_m$ must be changed from $\al$ in  
$\pi'_{r+2m}$ because the dependence on  $\nu_m$  is via
the product $\nu_m\al_m=\nu\al$. 

\begin{theorem} \label{prop:catal}
(i) For  $m\ge 0$ and $0\le r\le $ Ceiling $[1/\nu]$ \,, one has:
{\small
\begin{align}
\p'_{r,m}=\pi'_{r,m}\, 
p^{r+m}q^m\binom{r+2m}{m}\frac{r+1}{r+m+1}\,.
\end{align}
}

(ii) Let $\p_{r,m}^0\equal\lim_{\nu\to 0}\pi'_{r,m}$. Then
 $\p^0_{r,m}=\frac{\al^{r+2m}p^{r+m}q^m e^{-\al}}
{m! (r+m)!}\frac{r+1}{r+m+1}.$ I.e. we
count Catalan paths with the number of steps satisfying the 
$\al$-Poisson 
distribution. Furthermore, let  $p\!=\!\frac{1}{2}\!=\!q$ and
$\si_r^0\equal\sum_{m=0}^\infty \p^0_{r,m}$. Then 
\, $\si_r^0=\frac{r+1}{\al/2}e^{-\al}I_{r+1}(\al)$ 
for the {\sf hyperbolic Bessel function}
$I_u(z)\!=\!\sum_{m=0}^\infty 
\frac{(z/2)^{2m+u}}{m! \Ga(m+u+1)}$. Note that
$I_0(z)\!=\!\sum_{m=0}^\infty 
\frac{(z/2)^{2m}}{(m!)^2}$, $I_1(z)=dI_0(z)/dz$.

(iii) Under $p\!=\!1/2\!=\!q$ and  for $\nu=0$ as in (ii),  
the probability $\si^0_\ast=\sum_{r=0}^\infty \si_r^0$
to obtain a Catalan path with arbitrary  $r\ge 0$ and $m\ge 0$ 
equals $e^{-\al}\bigl(I_0(\al)+I_1(\al)\bigr)$. 
Also, the average value
$T=\sum_{r=0}^\infty (r\!+\!1)\,\si_r^0$ of the parameter
\,$(r\!+\!1)$ 
over such Catalan paths is $T=1$.

(iv) Setting $\nu\!=\!0$,  
let us consider {\sf all} paths, i.e. we disregard the
Catalan condition and $r$ can be negative, and let $p_r(m)$ be
the corresponding probability with $m$ steps and the output $r$.
 Then 
\begin{align}\label{prmall}
p_r(m)&=
p^{r+m}q^m\,\binom{r\!+\!2m}{m}
\frac{e^{-\al}}{(r\!+\!2m)!} \text{\ \, for } r\ge -m
\text{\,,\ and\ }\notag\\
&
P_r=\sum_{m=0}^\infty p_r(m)=
e^{-\al} \bigl(\frac{p}{q}\bigr)^{r/2}\, I_r\bigl(2\al(pq)^{1/2}\bigr).
\end{align}
In particular, $\sum_{r=-\infty}^\infty P_r=1$.
\end{theorem}
{\it Proof}. {\sf (i,ii)\,} The number of sequences $\{a_i=\pm 1,
1\le i\le n\}$ such that $\sum_{i=1}^k a_i\ge 0$ for 
any $1\le k\le n$
and $\sum_{i=1}^{n} a_i =r\ge 0$ is as follows: 
$C_{n}^{0\to r}=\binom{n}{(n-r)/2}\frac{r+1}{(n+r)/2+1}$,
where $n\!-\!r$ must be even. This is a standard formula
in the theory of Catalan numbers.
For $n=r+2m$, we obtain 
$C_{r+2m}^{0\to r}=\binom{r+2m}{m}\frac{r+1}{r+m+1}$. 
Then we multiply the latter expression by
$p^{r+m}q^m$ and by $\pi'_{r,m}$. This proves $(i)$; the justification
of $(ii)$ is straightforward.
\vskip 0.2cm

{\sf (iv)\,}  
To allow any sequences of
$\pm 1$, we must omit the terms 
$\frac{r+1}{m+r+1}$. Then $\sum_{r=\infty}^\infty P_r=1$,
which is the classical 
identity 
\begin{align}\label{bessel-id}
e^{\frac{\al}{2}
(t+\frac{1}{t})}=
\sum_{r=-\infty}^\infty I_r(\al)t^r, \text{ where }
I_r(\al)=I_{-r}(\al).
\end{align}

Indeed, $p_r(m)+p_{-r}(m)$ for  $r>0$ equals
{\small
$$(p^{r+m}q^m\!+\!q^{r+m}p^m)\binom{r\!+\!2m}{m}
\frac{e^{-\al}}{(r\!+\!2m)!}\!=\!
(pq)^{\frac{r+2m}{2}}\Bigl(\bigl(\frac{p}{q}\bigr)^{\frac{r}{2}}\!+\!
\bigl(\frac{q}{p}\bigr)^{\frac{r}{2}}\Bigr)\frac{e^{-\al}}{m!(r\!+\!m)!}.
$$
}
\!\!The sum of these terms and the one for $r=0$ is 
{\small
\begin{align*}
&P_r=e^{-\al}I_0(2\al(pq)^{1/2})+
e^{-\al}\sum_{r=1}^\infty \Bigl(\bigl(\frac{p}{q}\bigr)^{r/2}+
\bigl(\frac{q}{p}\bigr)^{r/2}\Bigr)I_r(2\al(pq)^{1/2})\\
&=
\exp\left(-\al+\al(pq)^{1/2}\Bigl(\bigl(\frac{p}{q}\bigr)^{1/2}+
\bigl(\frac{q}{p}\bigr)^{1/2}\Bigr)\right)=1
\text{ due to } p+q=1.
\end{align*}
}
\!\! Since we know that the sum must be $1$,  {\em this proves 
formula (\ref{bessel-id})}. Indeed, we can make $t=(p/q)^{1/2}$
arbitrary here.

{\sf (iii)\,}
Differentiating (\ref{bessel-id}) with respect to $t$, 
one obtains that
$
e^{\frac{\al}{2}(t+\frac{1}{t})-\al}=
\sum_{r=0}^\infty \p_r^0\, \frac{t^{r+1}-t^{-r-1}}{t-t^{-1}}$.
By making here $t\to 1$, we arrive at the formula $T=1$. Indeed,
$\frac{t^{r+1}-t^{-r-1}}{t-t^{-1}}=t^r+t^{r-2}+\cdots+t^{-r}$.
To obtain $\si_\ast^0$,
observe that the latter sum contains
$1$ for even $r$ and $t$ for odd $r$; then use
(\ref{bessel-id}) for $r=0,1$. 
\vskip 0.2cm

Actually, we do not need (\ref{bessel-id}) for $(iii)$.
The total number
of {\sf all\,} Catalan paths with $n$ steps is
$\binom{n}{n/2}$ for even $n$ and $\binom{n}{(n-1)/2}$
for odd $n$. It is $C_n^{0\to \ast}=
\sum_{r=0}^n C_{n}^{0\to r}$, where $n\!-\!r$ is even, in the
notation above. Thus,
the total probability to obtain a Catalan path is 
$$\sum_{r=0}^\infty \si_r^0=\sum_{n=0}^\infty
e^{-\al}\frac{\al^n}{n!}\binom{n}{\text{Floor}[n/2]}= 
e^{-\al}\bigl(I_0(\al)+I_1(\al)\bigr).
$$
Then, $2^n-\binom{n}{\text{Floor}[n/2]}=\sum_{r=0}^\infty
r\,C_{n}^{0\to r}$ for even $n\!-\!r$, since the ideal roulette is 
a $0$-sum
game, and $2^n=\sum_{r=0}^\infty
(r+1)\,C_{n}^{0\to r}$  using the formula above
for $C_n^{0\to \ast}$. This results in $T=1$.\sq

\begin{proposition}
Let $p\!=\!1/2\!=\!q$ and $\nu=0$ as in Part (iii) of
Problem \ref{prop:catal}. Now the starting level will
be $k\ge 0$ and the corresponding probability of Catalan paths will
be $
\si_\ast^0\lan k\ran=e^{-\al}\sum_{n=0}^\infty
\frac{\al^n}{n!}C_n^{k\to\ast}$. Then
{\small
$$
C_n^{k\to\ast}=
\sum_{n=0}^\infty \left(\binom{n}{\text{\rm\footnotesize 
Floor}[\frac{n-k}{2}]}\!+\!
\binom{n}{\text{\rm\footnotesize Floor}[\frac{n-k}{2}]\!+\!1}\!+\!\cdots\!+\!
\binom{n}{\text{\rm\footnotesize Floor}[\frac{n+k}{2}]}\right).
$$
}
\!\!\!Accordingly, 
 $\si_\ast^0\lan k\ran\!=\!
e^{-\al}\bigl(I_{-k}(\al)+I_{-k+1}+\cdots+I_{k+1}(\al)\bigr)$. 
\sq \end{proposition}

\comment{
{\small
$$
\si_r (\al)\sim 
\frac{1}{\sqrt{2\pi \al}}\left(1-\frac{4r^2-1}{8\al}+\ldots\right).
$$
} 
\!\!This  is due to the classical formula for 
$e^{-z}I_u(z)$ as $z\to \infty$.
For  $\al\sim 0$:\  
$\si_r(\al)\sim (\al/2)^r/r!$. 
}

\section{\bf Generating functions}
We will provide the generating function
for $\pi_r(N,L)$. We will show now the dependence of $\pi_r$
on $N,L$. Let 
$G(t,u)\equal\sum_{N=0}^\infty t^N \pi_r(N,L) t^N u^r$.
We fix $L$ here and below. 

\begin{theorem}
$G(t,u)=
\frac{1}{1-t}+\frac{\be\,(u-1)\,t}{(1-t)(1-(1-\be)t-\be u t^{L+1})}.$
\end{theorem}
{\it Proof.} First of all, let us calculate
$G^\circ_N=\sum_{r=0}^\infty  P(r,N,\be) u^r=
\sum_{r=0}^\infty  u^r \be^r
(1-\be)^{N-Lr-r} \binom{N-Lr}{r}$. This is the classical problem
about tiling the segment with $N$ boxes by $(L\!+\!1)$-minos, 
sequences
of $L\!+\!1$ consecutive
boxes,  and with $1$-minos. Its 
variant in a 2D square lattice with dominos 
and monominos 
(dimers and monomers)
is important in statistical physics. Though
there are no exact 2D formulas in the presence of $1$-minos. 
Here we count the tilings with the weights as above. For $u=1,\be=1/2,
L=1$:\, $G^\circ_N=f_{N+1}/2^N$ for the
Fibonacci numbers $f_N$. Generally: 
{\small
$$
G^\circ_N=(1-\be)G^\circ_{N-1}+\be u G^\circ_{N-L-1},\  
G^\circ_0=1,\, G^\circ_i=(1-\be)^i 
\text {\, for
\, } 1\le i\le L. 
$$
}
\!\! For instance, $G^\circ_{L+1}=(1-\be)^{L+1}+\be u.$
Using the standard facts in the theory of generating functions
or a straightforward  consideration:
$$
G^\circ(t,u)=\sum_{N=0}^\infty G^\circ_N t^N=
\frac{1}{1-(1-\be)t-\be u t^{L+1}}.
$$
Due to formula (\ref{Pform}),
$G_N\!=\!\sum_{r=0}^\infty  \pi_r(N,L) u^r$ 
satisfies the same recurrence as for $G_N^\circ$, but  with
different initial conditions. Namely,
$G(t,u)=(1+\be u t +\be u t^2+\ldots+ \be u t^{L})\,G^\circ(t,u)$. 
Finally, 
{\footnotesize
\begin{align*}
&G(t,u)=\frac{1-t+\be u t (1-t^{L})}{(1-t)(1-(1-\be)t-\be u t^{L+1})}\\
=&\frac{1-(1-\be)t-\be u t^{L+1} +\be u t -\be t}
{(1-t)(1-(1-\be)t-\be u t^{L+1})}=\frac{1}{1-t}\biggl(1+
\frac{\be\,(u-1)\, t }{(1-(1-\be)t-\be u t^{L+1})}\biggr).\text{\sq}
\end{align*} 
}

For $u=0$: $\sum_{N=0}^\infty \pi_0(N,L)=\frac{1}{1-(1-\be)t}$,
which we know without any calculations. For $u=1$:
$
\sum_{r=0}^\infty \pi_r(N,L)t^N=\frac{1}{1-t},
$
which gives a combinatorial proof of the identities
$\sum_{r=0}^\infty \pi_r (N,L)=1$ for any $N,L$. 
\vskip 0.2cm

{\sf Explicit formulas.}
The theorem readily gives that
{\small
\begin{align}\label{piu}
\pi_r(N,L)=\frac{\,d^r}{du^r}\,\biggl(\frac{1}{1-t}\Bigl(1+
\frac{\be\,(u-1)\, t }{1-(1-\be)t-\be u t^{L+1}}\Bigr)\biggr)
(u\mapsto 0).
\end{align}
}
\!\!Performing the differentiation, we obtain 
the following ``telescopic-type" presentation of $\pi_r$.

\begin{corollary}
Let $\Pi_r(L)=\frac{\be^r t^{(L+1)r-L}}{(1-t)(1-(1-\be)t)^r}$
for $r\ge 1$ and $\Pi_0(L)=\frac{1}{1-t}$.
Then $\sum_{N=0}^\infty \pi_r(N,L) t^N=\Pi_r(L)-\Pi_{r+1}(L)$ for
\ $r\ge 1$. This immediately gives that 
$\sum_{r=0}^\infty \pi_r(N,L)=1$ for any $N,L$. \sq
\end{corollary}

One can use this corollary 
to make the formulas for $\pi_r$ quite explicit: directly
expressed in terms of the binomial coefficients. The sums
there can be calculated using the standard combinatorial
identities. 

\begin{corollary}\label{cor:pir}
Provided that $N\!-\!L(r\!-\!1)-r\ge 0$,
{\small
$$
\pi_r(N,L)\!=\!\be^r\sum_{s=0}^m \binom{s\!+\!r\!-\!1}{r-1}
(1\!-\!\be)^{\max\{m-L,s\}}, \ m\!=\!N\!-\!L(r\!-\!1)-r,
$$
}
\!\!where there are $L+1$ powers of $(1-\be)$  
with integral positive 
coefficients for $m\ge L$. These powers are
$(1-\be)^{m-L+k}$ for $k=0,1,\ldots,L$ and the
coefficients 
depend only on $k,r$
unless $k=0$. For $m<L$, these 
powers are $(1-\be)^k$, where $0\le k\le m$ and all coefficients
depend only on $k,r$. \sq
\end{corollary}

\section{\bf Two processes} 
It is quite possible that several
strains (point processes) can be present simultaneously. 
They can be generally
with different immunity intervals $L$ 
and $\be$.  
\comment{For instance, when the Delta strain 
and the Omicron strain overlapped for some
times, the available vaccines were significantly more efficient
for the former than for the later. Statistically,
this alone makes the immunity durations different.  
The transmissibilities were
different too (greater for the Omicron).

Generally, vaccinations 
$(a)$ make the chances to be infected smaller, and  $(b)$ 
increase immunity intervals for those vaccinated if they are
infected. Correspondingly, they are reflected in $\be, L$
and $\al,\nu$ (for sufficiently large $N$).  As we discussed
above, the chances  to have no infection during
the whole cycle and the chances of being infected exactly 1 time
must be known for this. These chances depend on 
various factors, but the corresponding $L=\nu N$ can be 
potentially used
to evaluate the efficacy of different vaccines
in different countries. 
}

Let  $\be_1$ and $\be_2$ be the probabilities of being 
infected by strain 1 and strain 2 during 1 day, 
assuming that that the simultaneous infections 
by 1 and 2  are negligible. We set $\be_0=\be_1+\be_2$.
 The corresponding immunity intervals after the infections 
will be $L_1,L_2$. Let $\pi_{r_1,r_2}(N, L_1,L_2)$ be the
probability to have $r_1$ cases for strain 1 and $r_2$ 
for 2.  Accordingly, we need to calculate the generating
function $G=\sum_{N=0}^\infty G_N t^N$, where
$G_N=\sum_{r_1,r_2=0}^{\infty} u_1^{r_1} u_2^{r_2}\,
\pi_{r_1,r_2}(N, L_1, L_2)$.
\vskip 0.2cm

Similar to the above consideration, the basic combinatorial problem
is now to count the number of coverings of an $N$-segment by 
non-overlapping $(L_1\!+\!1)$-subsegments,
$(L_2\!+\!1)$-subsegments, and $1$-subsegments (monomers).
One of the subsegment can go through $N$, the endpoint
of the $N$-segment. Then $G_N$ satisfies the
recurrence relation $G_N=(1-\be_0) G_{N-1}+\be_1 u_1
G_{N-1-L_1}+\be_2 u_2 G_{N-1-L_2}$, and:
\begin{align*}
&G=\frac{1+u_1\be_1\,t\,(1+t+\ldots+t^{L_1-1})+
u_2\be_2\,t\,(1+t+\ldots+t^{L_2-1})}
{1-(1-\be_0) t-u_1\be_1 t^{L_1+1}-u_2\be_2 t^{L_2+1}}\\
&=
\frac{1-t+u_1\be_1\,t\,(1-t^{L_1})+ u_2\be_2\,t\,(1-t^{L_2})}
{(1-t)(1-(1-\be_0) t-u_1\be_1 t^{L_1+1}-u_2\be_2 t^{L_2+1})}\\
&=
\frac{1}{1-t}\biggl(1+
\frac{\be_1\,(u_1-1)t +\be_2\,(u_2-1)t}
{1-(1-\be_0) t-u_1\be_1 t^{L_1+1}-u_2\be_2 t^{L_2+1}}\biggr).
\end{align*}
\vskip 0.2cm

{\sf Several processes.}
The latter formula can be readily extended to any number
of simultaneous processes (strains). 
For $\be_0=\sum_i \be_i$ in the natural
notation:
{\small
\begin{align*}
&G=
\frac{1}{1-t}\biggl(1+
\frac{\sum_i\be_i\,(u_i-1)\,t }
{1-(1-\be_0)\, t-\sum_i u_i\be_i\, t^{L_i+1}}\biggr).
\end{align*}
}


Using that the dependence of $u_1,u_2$ is linear in the
numerator and denominator, it is not difficult to perform
the necessary $u$-differentiations and calculate the
generating functions with fixed  $r_1,r_2$.
For instance, let  $u_1=u=u_2$,  $\pi_r(N,L_1,L_2)$
be the probability that $r_1+r_2=r$, and $\p_r(t)=
\sum_{N=0}^\infty \pi_r(N,L_1,L_2)t^N$. Then

{\small 
$$\p_1(t)=
\frac{t (\be_1+\be_2) \left( 1-t (1- \be_1- \be_2)
-\be_1 t^{ L_1+1}-\be_2 
t^{ L_2+1}\right)}{(1-t) 
(1-t ( 1-\be_1-\be_2))^2}.
$$
}

\!\! This is for 1 infection by any strain (from two). 
When $L_1=L=L_2$
and $\be=\be_0=\be_1+\be_2$, we arrive at the case of one type of 
infection. For $r>1$, this is generally not true: even if $L_1=L_2$,
the order of the strains in their sequences matters.
An explicit combinatorial formula for
$N>L_2\ge L_1$ is as follows:
{\footnotesize
\begin{align*}
&\be_0^{-1}\pi_{1}(N, L_1\le L_2)=x^{N-L_2-1}
\biggl(x^0(N-L_2)(1-\be_1)+x^1\bigl(1-(N-L_2+1)\be_1\bigr)\\ &
+x^2\bigl(1-(N-L_2)\be_1\bigr)+\ldots+x^{L_2-L_1-1}
\bigl(1-(N-L_1-1)\be_1\bigr)\\
&+x^{L_2-L_1}+x^{L_2-L_1+1}+\ldots
+x^{L_2}\biggr), \text{\normalsize\, where\, 
$x\equal 1-\be_1-\be_2$.}
\end{align*}
}

It becomes somewhat simpler combinatorially in terms of $x,\be_2$:
{\footnotesize
\begin{align*}
&\be_0^{-1}\pi_{1}(N, L_1\le L_2)\,=\,x^{N-L_2-1}
\biggl(x^0(N-L_2)\be_2+x^1(N-L_2+1)\be_2\\ 
+&x^2(N-L_2+2)\be_2+\ldots+x^{L_2-L_1-1}
(N-L_1-1)\be_2+ x^{L_2-L_1}(N-L_1)\\
+&x^{L_2-L_1+1}+\ldots
+x^{L_2}\biggr), \text{\normalsize\, where\, 
$x\equal 1\!-\!\be_1\!-\!\be_2,\ N\!>\!L_2\!\ge\! L_1$.}
\end{align*}
}

\noindent
There are $L_2\!+\!1$ powers of $x$ here; the terms
{\small $x^1\bigl(1-(N-L_2+1)\be_1\bigr)$,
$x^1(N-L_2+1)\be_2$
} are present only if 
$L_2\!>\!L_1\!+\!1$.
For $N\!\le\! L_2$, the number of terms is $N$: they are exactly the
{\sf top\,} $N$ terms in the formulas above. The sums of the binomial
coefficients in this formula can be readily calculated,
which is useful for obtaining
the limits as $N\to \infty$, when
$N\be_i\to \al_i$, $L_i/N\to \nu_i$ for
$i=1,2$.  We note that by setting $\be_1=0, \be=\be_2, L=L_2$ in
the 1st formula, we obtain the 2nd where $\be_2=0, \be= \be_1,
L=L_1$. This is our formula for $\pi_1(N,L)$. 

Similarly, one calculates $\pi_{\lan 1 \ran}=\pi_{1,0}$ and 
$\pi_{\lan 2\ran}=\pi_{0,1}$, which
are the coefficients of $u_1$ and $u_2$ of $G$;
$\pi_{\lan i \ran}$ depends only on 
$L_i$. For $N>L_i$:
{\small
\begin{align*}
&\pi_{\lan i \ran }(N, L_i)\,=\,x^{N-L_i-1}\be_i
\Bigl(x^0(N-L_i)+x^1+\ldots+x^{L_i}\Bigr),\ i=1,2.
\end{align*}
}
\!\!where $x\equal 1\!-\!\be_1\!-\!\be_2.$ Accordingly, the 
top $N$ terms must be taken if $N\le L_i$.   Obviously, 
$\pi_1(N, L_1\le L_2)=
\pi_{\lan 1\ran }(N,L_1)+\pi_{\lan 2\ran}(N,L_2)$.  
In the limit $\be_i N\to \al_i$ and
$L_i/N\to \nu_i$, we obtain for $\al_0=\al_1+\al_2$:
{\small
\begin{align*}
&\pi'_{\lan i\ran}=\lim_{N\to\infty}\pi_{\lan i\ran}=
e^{\al_0(\nu_i-1)}\al_i
\Bigl(1-\nu_i+\frac{1-e^{-\al_0\nu_i}}{\al_0}\Bigr),\ i=1,2.
\end{align*}
}

Similar to Corollary \ref{cor:der} for $r=1$:
$$\frac{d \pi'_{\lan i\ran}}{d\, \nu_i}=
\al_i\al_0 e^{\al_0(\nu_i-1)}(1-\nu_i),
$$
We see that $\pi'_{\lan i\ran}$ increases
in terms of the corresponding $\nu_i$ for 
$0\le \nu_i<1$ and  fixed $\al_i>0$. This could be expected:
the greater $\nu_i$ the smaller the
total number of the corresponding 
reinfections. This is because  $\pi'_0$ does not depend on 
the immunity intervals.

\section{\bf Generalizations}
Let us now provide the generalizations of Theorems
\ref{thm:Fdif}, \ref{thm:Fdif-lim} to the case of any 
$L_i$ for $i=1,2,\ldots, s$ and the corresponding $\be_i$.
We set $\be_0=\be_1+\ldots \be_s$, $\be_i=\al_i/N$ and
$\al_0=\al_1+\ldots \al_s$. Let $\pi_{r_1,\ldots,r_s}$ 
be the probability as above,
where  $r_i$ is the number of segments of the 
length $L_i\!+\!1$ in $[1,N]$ corresponding to $\be_i$. 
One has: $\pi_{0,\ldots,0}=
(1-\be_0)^N$. We will assume below that $r\equal
r_1+\ldots+r_s$ is positive. Note that in the formulas below,
we set $\be_i^{r_i}=1=\al_i^{r_i}$ for $r_i=0$ even when
$\be_i=0=\al_i$. 

\begin{theorem}\label{thm:Fdif-new}
Setting $\si=N-\sum_{i=1}^s r_i L_i$,\, 
for $\si\ge r$:
\begin{align}\label{pipims}
\pi_{r_1,\ldots,r_s}&=\frac{
\bigl(\be_1^{r_1}\cdots \be_s^{r_s}\bigr)
(1-\be_0)^{\si-r}\, \si !} {\bigl(r_1!\cdots r_s!\bigr)\,(\si-r)!}\\
+\frac{\be_1^{r_1}\cdots \be_s^{r_s}}
{r_1!\cdots r_s!}\,
&\frac{d^{r-1}}{dX^{r-1}}
\biggl(\sum_{k=1}^{s} r_k\,X^{\si}\frac{1-X^{L_k}}{1-X}\biggr)
\bigl(X\mapsto 1-\be_0\bigr).\notag
\end{align}
Provided that $\be_i>0$ for any $i$,
the nonzero $\pi_{r_1,\ldots,r_s}$ are exactly when
$\si-r\ge L_k$ for at least one $k$. Accordingly,
for $\si<r$:
{\small
\begin{align}
&\pi_{r_1,\ldots,r_s}=
\frac{\be_1^{r_1}\cdots \be_s^{r_s}}
{r_1!\cdots r_s!}\,
\frac{d^{r-1}}{dX^{r-1}}
\biggl(\sum_{k=1}^{s} r_k\,X^{r-1}\frac{1-X^{\si+L_k-r+1}}{1-X}\biggr)
\bigl(X\mapsto 1-\be_0\bigr),\notag
\end{align}
}
\!\!where the summation is only over $k$ such that $L_k\ge r-\si$.
\sq
\end{theorem}

\begin{theorem}\label{thm:lim-new}
 We assume
that $L_k/N\to \nu_k$ for $1\le k\le s$ as $N\to \infty$.
Let $\pi'_{r_1,\ldots,r_s}\!=\lim_{N\to \infty} 
\pi_{r_1,\ldots,r_s}$,\, $\mu=1\!-\!
(\nu_1 r_1+\ldots+\nu_s r_s)$ and $\mu_k=\mu+\nu_k$. 
Then $\pi'_{0,\ldots,0}=e^{-\al_0}$. For $r>0$ and $\mu\ge 0$:
\begin{align}\label{piprims}
&\pi'_{r_1,\ldots,r_s}=\frac{
\mu^r \al_1^{r_1}\cdots \al_s^{r_s}}{r_1!\cdots r_s!}
e^{-\al_0 \mu}\\
+\,\frac{ \al_1^{r_1}\cdots \al_s^{r_s}\,(r-1)!}
{\al_0^{r}\,r_1!\cdots r_s!}&\sum_{k=1}^s\sum_{i=0}^{r-1}
\frac{r_k}{i!}\biggl(e^{-\al_0 \mu}(\al_0\mu)^i-
e^{-\al_0 \mu_k}(\al_0\mu_k)^i\biggr)\notag \\ 
=\,\frac{ \al_1^{r_1}\cdots \al_s^{r_s}\,(r-1)!}
{\al_0^{r}\,r_1!\cdots r_s!}&\sum_{k=1}^s
\biggl(\sum_{i=0}^{r}
\frac{r_k}{i!}e^{-\al_0 \mu}(\al_0\mu)^i-\!
\sum_{i=0}^{r-1}\frac{r_k}{i!}e^{-\al_0 \mu_k}(\al_0\mu_k)^i
\biggr).\notag
\end{align}
Additionally, the nonzero values occur
if $-\nu_k\!<\!\mu\!\le\! 0$ for at least one $k$:
{\small
\begin{align}\label{piprimse}
&\pi'_{r_1,\ldots,r_s}
=\frac{ \al_1^{r_1}\cdots \al_s^{r_s}\,(r-1)!}
{\al_0^{r}\,r_1!\cdots r_s!}\sum_{k=1}^s
\biggl(r_k-\!
\sum_{i=0}^{r-1}\frac{r_k}{i!}e^{-\al_0 \mu_k}(\al_0\mu_k)^i
\biggr),
\end{align}
}
\!\!\! where the summation is over $k$ such that  $\mu_k=\mu+\nu_k>0$.
If $\mu=0$ here, then (\ref{piprimse}) coincides with the
formula above; recall that $0^0=1$. 
\end{theorem}
{\it Proof.} For the first term in (\ref{pipims}):
{\small
$$
\frac{
\bigl(\be_1^{r_1}\cdots \be_s^{r_s}\bigr)
(1-\be_0)^{\si-r}\, \si !} {\bigl(r_1!\cdots r_s!\bigr)\,(\si-r)!}
\to \frac{
\mu^r \al_1^{r_1}\cdots \al_s^{r_s}}{r_1!\cdots r_s!}
e^{-\al_0 \mu}.
$$
}

Following the proof of Theorem \ref{thm:Fdif-lim}, the remaining
summation in (\ref{pipims})
becomes the limit of 
{\small
$$\sum_{k=1}^s \sum_{i=0}^{r-1}
\frac{r_k \be_1^{r_1}\cdots\be_s^{r_s}}{r_1!\cdots r_s!}
\binom{r-1}{i}\Bigl(\si^i-(\si+L_k)^i\Bigr)
\frac{(r-1-i)!}{\be_0^{r-i}}\,.
$$
}
\!\!Then we perform the substitutions $\be_i=\al_i/N$ and
$L_i=\nu_i N$, $\be_0^i\si^i=(\al_0\mu)^i$,
$\be_0^i(\si+L_k)^i=(\al_0\mu_k)^i$ and tend $N$ to $\infty$. 
\vskip 0.2cm

When $\mu<0$, the above sum becomes:
{\small
$$\sum_{k=1}^s \sum_{i=0}^{r-1}
\frac{r_k \be_1^{r_1}\cdots\be_s^{r_s}}{r_1!\cdots r_s!}
\binom{r-1}{i}\Bigl((r\!-\!1)\cdots(r\!-\!i)-(\si+L_k)^i\Bigr)
\frac{(r-1-i)!}{\be_0^{r-i}}\,,
$$
}
\!\!\!where $(r\!-\!1)\cdots(r\!-\!i)\be_0^i$ vanishes in the limit
unless $i=0$. The summation is over $k$
such that $\si+L_k\ge r$, which gives  $\mu_k>0$. 
\sq

Assuming that $0<\nu_1\le\nu_2\le \ldots \le\nu_s$, we can define
$\kappa\lan \mu\ran=\min\{k \mid \mu+\nu_k>0\}$ for $-\nu_s<\mu\le 0$;
the corresponding $\{r_k\}$ are called {\sf extreme}.
The condition   $\mu_k=\mu+\nu_k>0$ 
in (\ref{piprimse}) becomes $k\ge \kappa\lan \mu \ran$ and this
formula reads for  $-\nu_s<\mu\le 0$:
{\small
\begin{align}\label{piprimsee}
&\pi'_{r_1,\ldots,r_s}
=\frac{ \al_1^{r_1}\cdots \al_s^{r_s}\,(r-1)!}
{\al_0^{r}\,r_1!\cdots r_s!}\sum_{k=\kappa\lan \mu \ran}^s
\biggl(r_k-\!
\sum_{i=0}^{r-1}\frac{r_k}{i!}e^{-\al_0 \mu_k}
(\al_0\mu_k)^i
\biggr). 
\end{align}
}
\!\!The extreme configurations have little to do with epidemics,
but can have potential applications in physics and networks. 
For instance, let us estimate  $\pi'_{r_1,\ldots,r_s}$ at the
{\sf resonances}: when $\mu= 0$. Then 
$\kappa=1$ and  the
Taylor formulas gives that for some $0< \ep_k <\nu_k$:
$$
\frac{ \al_1^{r_1}\cdots \al_s^{r_s}}
{r_1!\cdots r_s!}e^{-\al_0 \nu_s}<\, \pi'_{r_1,\ldots,r_s}
\!\!=\! \frac{ \al_1^{r_1}\cdots \al_s^{r_s}}
{r_1!\cdots r_s!}\sum_{k=1}^s \frac{r_k}{r}
e^{-\al_0 \ep_k}< \frac{ \al_1^{r_1}\cdots \al_s^{r_s}}
{r_1!\cdots r_s!}.
$$

The dependence on $\mu$ is of interest. It takes values
in the set $\{\mu\}=\{1-\sum_{k=1}^s r_k\nu_k >-\nu_s \mid r_k\ge 0\}$;
we can define $\pi'(\mu)$ as the sum of $\pi'_{r_1,\ldots, r_s}$
over the configurations
with $1-\sum_{k=1}^s r_k\nu_k=\mu$. The sequence $\{r_k\}$
is uniquely determined
by $\mu$ if $\nu_k$ are in a general position.
For rational $\nu_k=\ell_k/N$:
$\{\mu\}=\{N-\sum_{k=1}^s r_k \ell_k >-\ell_s \mid r_k\ge 0\}/N$.
 The extreme $\mu$
are then for $\{r_k\}$ such that $N+\ell_s>\sum_{k=1}^s r_k \ell_k\ge N$,
which is linked to the Frobenius coin problem.
\comment{
For instance, let $s=2$ and $\nu_1\ll \nu_2$. 
The most likely extreme configuration will be with 
$r_2\approx 1/\nu_2$ and $0<r_1 <\frac{\nu_2}{\nu_1}$. 
Let us take
$\al_1=\al=\al_2$ and assume that $\ep_k\approx \nu_k$.
Qualitatively, we need to maximize 
$\al^{1/\nu_2}e^{-2\al\nu_1}$; its 
point of maximum is
$\al=\frac{1}{2\nu_1\nu_2}$. 
This is very approximate.
However, this calculation demonstrates that the problem
of maximization of the extreme probabilities 
makes sense. 
}

\vskip 0.2cm
{\sf Extreme configurations.}
There is a natural way to obtain a distribution of
probabilities {\sf only on the extreme configurations}.
We set $\rho_i=\al_i/\al_0$ and  consider
$\tilde{\pi}'_{r_1,\ldots, r_s}\equal \lim_{\al_0\to \infty}
\pi'_{r_1,\ldots, r_s}$. For $r>0$:
{\small
\begin{align} \label{piextr}
&\tilde{\pi}'_{r_1,\ldots,r_s}=
\frac{\rho_1^{r_1}\cdots \rho_s^{r_s} r!}
{r_1!\cdots r_s!}
\sum_{k=\kappa\lan \mu\ran}^s \frac{r_k}{r},\ \  
\mu=1-\sum_{k=1}^s r_k \nu_k,
\end{align}
}
\!\!for extreme $\{r_k\}$, and $0$ otherwise. For finite $\al_0$,
all $\mu>0$ occur; generally, the ``defects" (gaps)
 play a very significant role.  
\vskip 0.2cm

The combinatorial
counterpart of (\ref{piextr}) is straightforward. We assume that 
$\be_0$ is significantly larger than $1/N$ and that 
the powers $(1-\be_0)^{\nu_k N}$ can be disregarded.
This means that we cover $[1,N]$ by
$L_k$-segments {\sf without gaps} between them (the last segment
can go beyond $N$). The lengths
of the segments are here $\ell_k\ge 1$, not $(L_k\!+\!1)$ as in 
Theorem \ref{thm:Fdif-new}. Then $\be_0=1$ and 
we set $\rho_k=\be_k$. Let $0< \ell_1\le \ell_2\le \ldots
\le \ell_s$; as above, $\si=N-\sum_{k=1}^s r_k \ell_k,\ 
\si_k=\si+\ell_k$. 
Then, provided that $-\ell_s< \si \le 0$,
$$
\tilde{\pi}_{r_1,\ldots,r_s}=
 \frac{\rho_1^{r_1}\cdots \rho_s^{r_s} r!}
{r_1!\cdots r_s!}
\sum_{k=\text{\tiny{$\kapp$}}\lan \si\ran}^s \frac{r_k}{r},\ 
\text{\small{$\kapp$}}\lan \si \ran=\min\{k \mid \si_k>0\},
$$
which 
coincides with (\ref{piextr}) in the limit $\ell_i/N\to \nu_i$.
We obtain that
$\sum_{\{r_k\}}\tilde{\pi}_{r_1,\ldots,r_s}=1$ and deduce from it 
that $\sum_{\{r_k\}}\tilde{\pi}'_{r_1,\ldots,r_s}=1$.
\vskip 0.2cm
 
{\sf Correlation functions.} 
As above, we will cover  $[1,N]$ by $r=r_1+\ldots+r_s$ segments of lengths
$\ell_1,\ldots, \ell_s$ (with probabilities $\rho_1,\ldots,\rho_s$).
The corresponding $m$-point function
is defined for  $m\ge 1,\, 1\le k_i \le s$, and  
$1\le N_1\le N_2\le \ldots\le N_m$.
It is the probability 
$\p^{k_1,\ldots, k_m}_{N_1,\ldots,N_m}(r_1,\cdots, r_s)$
of such coverings subject to the following:
{\small
$$ N_1\in X[k_1] \text{\,\, and\, } 
\ N_i+\ell_{k_1}+\ldots \ell_{k_{i-1}} \in X[k_i]
\text{\,\, for any\, } 1<i \le m,$$
}
\!\!where $X[k]$ is the corresponding consecutive segment of 
the configuration; $k$ denotes its length: $|\,\!X[k]\,\!|=\ell_k$.
For instance, the $1$-point function is when 
$X[k_1]$ covers $N_1$ such that $1\le N_1\le.$

Let $r_k^{(i)}\ge 0$ for $i=1,2,\ldots, m+1$ be the number of 
segments of types $\,k=1,2,\ldots,s\,$  strictly between
 $X[k_{i-1}]$ and $X[k_{i}]$,
before $X[k_1]$ for $i=1$, and (strictly) after
$X_{k_{m}}$. Thus, 
 $r_k\!-\!\sum_{i=1}^{m+1} r_k^{(i)}$ is
the number of $k_i=k$
in the sequence $\{k_i, 1\!\le\! i\!\le\! m\}$; this relation will be
imposed below.
We set
$\varsigma^{(i)}\equal \sum_{j=1}^i (r_1^{(j)}\ell_1+
\ldots+r_s^{(j)}\ell_s)$ for $1\le i\le m\!+\!1$.
For instance, $\varsigma^{(m+1)}=\sum_{i=1}^m r_i\ell_i-
(\ell_{k_1}+\ldots+\ell_{k_m})$.
We obtain the following formula
in terms of the multinomial coefficients $\binom{a+b+c+\ldots}
{a,\,b,\,c,\ \ldots}$:
{\footnotesize
\begin{align}\label{partext}
&\ \ \ \ \ \ \text{\normalsize
$\p^{k_1,\ldots, k_m}_{N_1,\ldots,N_m}(r_1,\ldots, r_s)$}\ 
=\ \sum \rho_1^{r_1}\cdots \rho_s^{r_s}
\binom{r_1^{(1)}+\ldots+r_s^{(1)}}
{r_1^{(1)},\ \ldots\,,\, r_s^{(1)}}\cdots\\
&\cdots
\binom{r_1^{(m)}+\ldots+r_s^{(m)}}
{r_1^{(m)},\ \ldots\,,\, r_s^{(m)}}
\binom{r_1^{(m+1)}+\ldots+r_k^{(m+1)}+\ldots+r_s^{(m+1)}\!-\!1}
{r_1^{(m+1)},\,\ldots\,,\,r_k^{(m+1)}-1\,,\,\ldots\,,\,r_s^{(m+1)}},\notag\\
&\text{\normalsize where} \ \text{\small the summation is over  
\,$1\le k\le s$\, and\, $\{r_k^{(i)}\ge 0\}$\, such that}
\notag\\
&N_1+\ell_{k_1}>\varsigma^{(1)}\!\ge\! N_1,
N_2+\ell_{k_2}>\varsigma^{(2)} \!\ge\! N_2,\ \ldots,\ 
N_m+\ell_{k_m}>\varsigma^{(m)}\ge N_m, \notag\\
&\text{\small
and 
$N\!+\!\ell_k>\varsigma^{(m+1)}\!+\!\ell_{k_1}\!+\ldots+\!\ell_{k_m}
\ge N$, equivalently, $ N\in  X[k]$.} \notag
\end{align}
}
\!\!\!The last term in the product occurs only if 
$r_1^{(m+1)}+\ldots+r_s^{(m+1)}\ge 1,$ i.e. when
the last segment is not fixed in the $m$-point function. 

The $m$-point function serving Theorem \ref{thm:Fdif-new} is when
$\rho_k=\be_k,\, \ell_k= L_k+1$, and we add
$\ell_0=1,\, \rho_0=1-\be_0$.  It is 
$\p^{k_1,\ldots, k_m}_{N_1,\ldots,N_m}(r_0,r_1,\ldots, r_{s})$,
where now $k_i\ge 0$. We can omit $r_0$ here, the number
of positions not covered by the segments, and fix only the 
corresponding numbers of segments: $r_1,\ldots, r_s$.
We will not discuss the limit as $N\to\infty$, providing only
the following simple example.

Let $s=1$, so we have only $1$ length $L$. Accordingly $\be=\al/N$,
and we use below the formulas for $\pi_r'$ from
Theorem \ref{thm:Fdif-lim}
in terms of $\be, N, L$:
$\al\mapsto \be N, \nu\mapsto L/N$. We will write 
$\pi'_r(N; L,\be)$.
 
Then for $m=1$, the limit $N\to \infty$ 
of $\p^{0}{N_1}(r)$, which is
when the position $N_1\le N$ is not covered by a segment, is  
\begin{align*}
&\p^{\prime\,0}_{N_1}(r)=
\sum_{i_1,i_2} e^{(i_1L-N_1)\be} \pi'_{i_2}(N-N_1;L,\be) 
\for \\
& 0\le i_i <\frac{N_1+1}{L+1},\, 
0\le i_2 \le \frac{N_1+1}{L+1},\, i_1+i_2=r.
\end{align*}
Here we do not need the limiting procedure; this is a 
direct calculation with probabilities.
Also, $\p^{\prime\,1}_{N_1}(r)$, the 1-point function
when $N_1$ is covered by an $(L+1)$-segment,
is $\pi'_r(N; L,\be)-\p^{\prime\,0}_{N_1}(r)$. 

\comment{
We note that $\pi_{r_1,\ldots,r_s}$ and 
$\pi'_{r_1,\ldots,r_s}$ are nonzero if and only if
$\si+L_k=N-\sum_{i=1}^s r_i L_i+L_k\ge r$ for at least one $k$.
Accordingly, $\mu_k\ge r$ for at least one $k$ is necessary and
sufficient for 
$\pi'_{r_1,\ldots,r_s}\neq 0$. The inequalities $\si\ge r$ and $\mu\ge r$
are more restrictive. The formulas for $\{r_i\}$ beyond the
latter inequalities can be readily calculated; we will omit them. 
\vskip 0.2cm 
}

\section{\bf Some perspectives}
The consideration of {\sf ensembles} of segments
of various lengths $(L_k+1)$ 
links our paper to {\sf stochastic processes}, namely
to {\sf Whittaker-type processes}, those based on 
the distances between {\sf neighboring} particles.
See \cite{BC}. 
In our approach, $L_k$ with the corresponding multiplicities $r_k$
are some substitutes
for these distances with one reservation: recall that 
we allow ``defects", the gaps between our segments.  
Thus, $\{r_k,L_k\}$ basically give these distances.
\vskip 0.2cm

Here the {\sf Mat\'{e}rn II statistics} is employed. Namely, 
for consecutive time moments $i=1,2,\ldots, N$,  
the segment $[i,i+L_k]$ can be created at any such $i$ with the
probability $\be_k$. This is unless $i$  belongs to the previously created 
segment, when the point is deleted (this operation of
{\sf thinning}). 
We allow
only finitely many  possible $L_k$: $1\le k\le s$. 
Generally, one can always expected some formulas 
in terms of binomial coefficients
for other kinds of statistics.
\vskip 0.2cm

There is an almost immediate link to the so-called 
{\sf interlacing sequences}
$x_1<y_1<x_2<y_2<\ldots<x_{n-1}<y_n<x_{n+1}$. They are
sequences of $n$ non-overlapping segments $[x_i,y_i]$ of 
lengths $y_i-x_i$ in $[x_1,x_{n+1}]$.  The 
corresponding {\sf transitional probabilities} are
associated with the $t$-residues of the function
$F(t)=\frac{(t-y_1)(t-y_2)\ldots 
(t-y_n)}{(t-x_1)(t-x_2)\ldots (t-x_{n+1})}$
at its poles: $x_1,\ldots, x_{n+1}.$ See \cite{Ke,BO,Ol}. In our approach,
these segments are $[x_i, y+i=x_i+L_{k_i}]$ where $L_{k_i}$ is the
length of the segment from $x_k$. Then we calculate the
generating function $G(t)$ and expand it; its coefficient
of $t^N u_1^{r_1}u_2^{r_2}\cdots$ is the probabilistic 
measure of the corresponding {\sf Young diagram} 
of order $\le N$; $r_k$ is the number of
rows of the length $L_k+1$. 
\vskip 0.2cm

For instance, one can take here
$L_i=i$ for $i=1,2,\ldots, n$. Then our approach becomes 
close to the theory of the transitional 
probabilities. The denominator of $G(t)$ is basically 
that in $F(t)$ with undetermined coefficients, which are 
our  $u_i\be_i$. The numerator 
of $G(t)$, a polynomial of degree $n$, 
  incorporates  the edge effects. It is not ``generic" in this
approach; this is different.  
The points $x_i, y_i$ do not appear; the aim is to obtain
the corresponding distribution of $\{r_k\}$ in terms 
of the probabilities
$\be_k$. The latter 
can depend on the corresponding $L_k$, and even  
on the whole configuration of
segments.
\vskip 0.2cm
  
The classical theory results in the distribution of probabilities
for Young diagrams
related to the celebrated hook-formula and {\sf Jack polynomials}. 
Our probabilities of Young diagrams
from (\ref{pipims}) and (\ref{piprims}) are of different
nature; we think that this approach is new. 
\vskip 0.2cm

A variant of our approach is for ``Poisson-Catalan processes,"
where $\pm L_i$ are considered as the jumps up and down of
the ``energy function". 
Theorem \ref{prop:catal} is stated only for one $L$ (for $s=1$).
The probabilities $\p'_{r,m}$ there 
can be extended to multi-dimensional (type $A_s$)
Catalan paths, where the
steps are at the 
points subject to the Poisson-type distribution
from Theorem \ref{thm:lim-new}. Generally, one needs here
formulas for the  number of
standard Young tableau for a given skew Young diagram, for
instance the  {\sf Naruse hook-length formula}. See e.g.
\cite{MPP}. Certain multi-dimension Bessel-type functions
occur here, which can be interesting to study.


\vskip 0.2cm

{\sf Acknowledgements.}
The author thanks very much Alexei Borodin for important
discussions, and Evgeny Feigin for his help. 

\comment{
{\sf Author's Affiliation.}
University of North Carolina at Chapel Hill, Mathematics Department,
Phillips Hall, Chapel Hill, NC 27599-3250.
}

\bibliographystyle{unsrt}

\begin{thebibliography} {ABCD}
\vbadness=10000

\bibitem [ADDP] {ADDP}
F.~Arruda, S.~Das, C.~Dias, D.~Pastore,
{\em Modelling and optimal control of multi strain 
epidemics, with application to COVID-19},
PLoS ONE 16: 9 (2021), e0257512; doi.org/10.1371/journal.pone.0257512.


\bibitem [BC] {BC}
A.~Borodin, I.~Corwin,
{\em Macdonald processes},
Probability Theory and Related Fields 158: 1 (2014), 225--400.

\bibitem [BO] {BO}
A.~Borodin, G.~Olshanski, 
{\em Point processes and the infinite symmetric group}, Math.
Research Lett. 5 (1998), 799--816. 

\bibitem  [Ch1] {ChB} 
I.~Cherednik,
{\em Momentum managing epidemic spread and Bessel \\
functions}, Chaos, Solitons \& Fractals {139} (2020);
doi.org/10.1016/\\ j.chaos.2020.110234. 

\bibitem  [Ch2] {ChM}
I.~Cherednik, {\em Modeling the waves of Covid-19},
Acta Biotheoretica 70, 8 (2022);
doi.org/10.1007/s10441-021-09428-w.


\bibitem [Gi] {Gi}
P.~Giaquinta, 
{\em Entropy and ordering of hard rods in one dimension},
Entropy 2008, 10, 248--260; 
doi.org/10.3390/e10030248.

\bibitem [Ke] {Ke}
S. V.~Kerov, 
{\em Anisotropic Young Diagrams and Jack Symmetric Functions}, 
Funktsional. Anal. i Prilozhen., 34: 1 (2000), 51--64; 
Funct. Anal. Appl., 34: 1 (2000), 41--51.

\bibitem [KD] {KD}
K.~Koufos, C.~Dettmann, 
{\em Moments of interference in
vehicular networks with hardcore headway distance}, IEEE
Transactions on Wireless Communications, 17 :12 (2018), 8330--8341;
doi.org/10.1109/TWC.2018.2876241.

\bibitem [MPP] {MPP}
A.~Morales, I.~Pak, G.~Panova,
{\em Hook formulas for skew shapes I. q-analogues and bijections},
Journal of Combinatorial Theory, Series A 154 (2018), pp 350--405.

\bibitem [Ol] {Ol}
G.~Olshansky,
{\em Random permutations and related topics},
The Oxford Handbook of Random Matrix Theory, 
Edited by Akemann, Baik, and Di Francesco (2015).




\end{thebibliography}

\end{document}